\documentclass[12pt]{amsart}

\usepackage{amssymb,latexsym}

\usepackage{xcolor}

\usepackage{enumerate}

\usepackage[T1]{fontenc}
\usepackage{hyperref}
\hypersetup{
    colorlinks=true,
    linkcolor=blue,
    filecolor=blue,      
    urlcolor=blue,
    citecolor=blue,
    pdftitle=,
    pdfpagemode=FullScreen,
    }

 \usepackage[french,english]{babel}


\newtheorem{thm}{Theorem}[section]

\newtheorem{lem}[thm]{Lemma}
\newtheorem{pro}[thm]{Proposition}
\theoremstyle{definition}

\newtheorem{rem}[thm]{Remark}

\numberwithin{equation}{section}

\renewcommand{\leq}{\leqslant}
\renewcommand{\ge}{\geqslant}
\renewcommand{\geq}{\geqslant}

\newcommand{\R}{\mathbb{R}}
\newcommand{\C}{\mathbb{C}}

\newcommand{\X}{\mathbb{X}}

\newcommand{\re}{\textup{Re}}
\newcommand{\im}{\textup{Im}}

\newcommand{\mc}{\mathcal}

\newcommand{\eps}{\varepsilon}

\newcommand{\newabstract}[1]{%
  \par\bigskip
  \csname otherlanguage*\endcsname{#1}%
  \csname captions#1\endcsname
  \item[\hskip\labelsep\scshape\abstractname.]
}

\begin{document}

\baselineskip=16.5pt

\title[]{Isolated zeros of the Riemann zeta function and zeros of linear combinations of Dirichlet $L$-functions}

\author{Youness Lamzouri}

\address{
Universit\'e de Lorraine, CNRS, IECL, 
F-54000 Nancy, France}

\email{youness.lamzouri@univ-lorraine.fr}


\begin{abstract} Let $0<\gamma_1\leq \gamma_2\leq \cdots$ be the ordinates of the non-trivial zeros of the Riemann zeta function in the upper half-plane, counted with multiplicity, and fix $c>0$. We say that a non-trivial zero $\rho_n=\beta_n+i\gamma_n$ is $c$-isolated if $\min(\gamma_{n+1}-\gamma_n, \gamma_n-\gamma_{n-1})>2\pi c/\log \gamma_n$.  
Note that any isolated zero is necessarily simple and on the critical line.  We prove unconditionally that there exists a positive constant $c>0$ such that at least $13/20$ of the zeta zeros are $c$-isolated. Moreover, we can show that this proportion is more than $1/2$ when $c=1/6$. 
Previously, it was not known unconditionally that there are infinitely many isolated zeros.  We also show that there are positive proportions of $k$ consecutive ``moderate'' gaps (gaps larger than a constant times the average local gap) between the ordinates of the zeta zeros, for $k=3, 4, 5$. 

As an application we combine our work on isolated zeros (which we extend to Dirichlet $L$-functions) with an unconditional refinement of the  Bombieri--Hejhal method to show unconditionally that at least $31\%$ of the zeros of a linear combination of two distinct primitive Dirichlet $L$-functions with the same conductor and parity lie on the critical line and have odd multiplicity. The previous unconditional record for this proportion was $7.38\times 10^{-9}$, which was obtained by Dousselin who refined and made explicit a result of Selberg. 
\end{abstract}


\maketitle


\section{Introduction}

The distribution of the non-trivial zeros of the Riemann zeta function $\zeta(s)$ is one of the central problems in analytic number theory. Although there are precise conjectures concerning their statistics, our knowledge remains limited, even under the assumption of the Riemann hypothesis (RH). In particular, the vertical distribution of these zeros is ultimately connected to the distribution of prime numbers in short intervals, and to the size of the class numbers of imaginary quadratic fields. 

Let $\rho_n=\beta_n+i\gamma_n$ denote the non-trivial zeros of the zeta function in the upper half-plane, ordered so that $0<\gamma_1\leq \gamma_2\leq \cdots$, where the zeros are counted with multiplicity. The classical Riemann--von Mangoldt formula gives  
$$ N(T):=\sum_{0<\gamma_n\leq T} 1= \frac{T}{2\pi}\log \left(\frac{T}{2\pi}\right)-\frac{T}{2\pi}+O(\log T).$$
Thus, the average spacing between consecutive ordinates of zeros $\gamma_{n+1}-\gamma_n$ is $2\pi/\log \gamma_n$. To measure how gaps deviate from their average, the following quantities have been extensively studied
$$\mu:=\liminf_{n\to \infty} \frac{\gamma_{n+1}-\gamma_n}{2\pi/\log \gamma_n}, \qquad \textup{and } \qquad \lambda:=\limsup_{n\to \infty} \frac{\gamma_{n+1}-\gamma_n}{2\pi/\log \gamma_n}.$$
One trivially has $\mu\leq 1\leq \lambda$. The best current result for small gaps is  $\mu\leq 1/2$ obtained independently by Bondarenko and Heap \cite{BonHea} and by Durkan and Pearce-Crump \cite{DurPea}, improving the bound $\mu<0.50895$ due to Inoue \cite{In}. For large gaps, the best results are obtained under the assumption of the Riemann hypothesis and state that $\lambda>3.63$ by Durkan and Pearce-Crump \cite{DurPea} improving the previous record $\lambda>3.18$ due to Bui and Milinovich \cite{BuMi}. On the other hand, the work of Conrey and Iwaniec \cite{CoIw} shows that small spacings between zeros of the zeta function are connected to the size of the class number of imaginary quadratic fields. 

It is widely believed that $\mu=0$ and $\lambda=\infty$. The first follows from Montgomery's pair correlation conjecture (PCC) which asserts that assuming RH we have 
\begin{equation}\label{Eq:PairCorrelationConj}
\frac{1}{N(T)}\#\left\{0<\gamma_n, \gamma_m\leq T : 0< \gamma_n-\gamma_m \leq \frac{2\pi \alpha}{\log T}\right\}= \mc{P}(\alpha)+o(1),  
\end{equation}
for any fixed $\alpha\geq 0$, where 
\begin{equation}\label{Eq:DefinitionPAlpha}
\mc{P}(\alpha):= \int_0^{\alpha} \left(1-\left(\frac{\sin \pi u}{\pi u}\right)^2\right) du.  
\end{equation}
The conjecture $\lambda=\infty$ is a consequence of the stronger GUE Hypothesis which asserts that the statistics of the ordinates of the non-trivial zeros of the zeta function, appropriately normalized, behave like the statistics of eigenvalues of large random matrices from the Gaussian Unitary Ensemble (GUE). Note that the pair correlation conjecture is the two-level case of the GUE Hypothesis.  

The pair correlation conjecture has strong implications for the distribution of prime numbers in short intervals. Indeed, Gallagher and Mueller \cite{GaMu} proved that assuming RH and PCC one has 
$$ \int_1^x \left(\psi\big((1+\delta) u\big)- \psi(u)-\delta u\right)^2\frac{du}{u^2}\sim \delta \log\left(\frac{1}{\delta}\right)\log\left(x\sqrt{\delta}\right),
$$
as $x\to \infty$, uniformly for $1/x\leq \delta \leq x^{-\varepsilon}$ for any fixed $0<\eps<1$, where $\psi(x)=\sum_{n\leq x} \Lambda(n)$, and $\Lambda(n)$ is the von Mangoldt function (which equals $\log p$
 if $n=p^k$ for a prime $p$ and positive integer $k$, and equals $0$ otherwise). 
Weaker forms of the pair correlation conjecture were used by several authors to obtain various applications (see for example \cite{BoHe}, \cite{Ci}, \cite{CiGoKi}). In particular,  Bombieri and Hejhal \cite{BoHe} used the following well-spacing hypothesis 
\begin{equation}\label{Eq:WellSpacingHyp}
\lim_{\eps\to 0+}\limsup_{T\to\infty}\frac{1}{N(T)}\#\left\{0<\gamma_n\leq T : 0\leq \gamma_{n+1}-\gamma_n \leq \frac{2\pi \varepsilon}{\log T}\right\}= 0,
\end{equation}
for general $L$-functions together with the generalized Riemann hypothesis GRH to prove that $100\%$ of the zeros of linear combinations of $L$-functions (satisfying certain natural conditions) lie on the critical line, confirming a conjecture of Montgomery (see Section \ref{Sec:LinearCombinations} below.)

\subsection{Isolated zeros and consecutive gaps between the zeta zeros} In \cite{GoSa} Gonek and Sahay asked the following basic question: can one prove the existence of infinitely many zeros of the zeta function that are ``isolated'', in the sense that both neighboring ordinates are away by at least a positive constant times the average gap? Prior to this work (and the independent work \cite{ARTZ}, see Remark \ref{Rem:Others} below), this was not known unconditionally. Note that isolated zeros are necessarily simple and on the critical line. Our first result answers Gonek and Sahay's question in the affirmative in a strong sense, showing unconditionally that more than $13/20$ of the zeros of the Riemann zeta function are isolated. 
\begin{thm}\label{Thm:Isolated0}
For $c>0$ let $B_0(T, c)$ be the number of  zeros $\rho_n$ with $n\geq 2$ such that $0<\gamma_n\leq T$ and  $\min\big(\gamma_{n+1}-\gamma_n, \gamma_{n}-\gamma_{n-1}\big)> 2\pi c/(\log T)$. 
There exists a positive constant $c>0$ such that 
\begin{equation}\label{Thm:IsolatedSmall}
\liminf_{T\to \infty} \frac{B_0(T, c)}{N(T)}\geq 0.6552..
\end{equation}
Moreover, we have 
\begin{equation}\label{Thm:Isolated}
\liminf_{T\to \infty} \frac{B_0(T, 1/6)}{N(T)}\geq 0.50875.
\end{equation}
\end{thm}
\begin{rem}\label{Rem:Optimization}
We note that we did not try to optimize the proportion of $c$-isolated zeros for small $c>0$ in \eqref{Thm:IsolatedSmall}, and we believe that some further small improvement is possible. There are two places where the method can be further pushed. First, one can optimize the kernel used in Proposition \ref{Pro:IsolatedSMALL} below. We used the Montgomery-Taylor kernel which is optimal for the pair correlation sum in Lemma \ref{Lem:ConstructionKernel} (see \cite{CCLM}). However, in our setting one needs to optimize the full term on the right hand side of \eqref{Eq:SMallGapsProof}, which contains a secondary term coming from off critical line zeros. To bound this secondary term we use a result of Conrey \cite{CoBLMS} from 1989, which uses his asymptotic formula for the mollified second moment of the zeta function together with Littlewood's lemma. Thus a second source of a possible improvement would be to use better mollifiers to optimize this secondary term. 
\end{rem}

Fix $c>0$. Gonek and Sahay \cite{GoSa} called a gap $\gamma_{n+1}-\gamma_n$ between consecutive ordinates ``moderate'' if $\gamma_{n+1}-\gamma_n\geq 2\pi c/\log\gamma_n$. Since this notion depends on $c$, we shall instead refer to it throughout as a $c$-moderate gap. Thus \eqref{Thm:Isolated} shows that at least $50.87\%$ of the consecutive gaps are $1/6$-moderate. Optimizing this result we can show that the same proportion of consecutive gaps is $c$-moderate where $c=0.1696.$
Gonek and Sahay asked the more general question: how frequently are there $k$ consecutive $c$-moderate gaps for a fixed small $c>0$? They showed assuming RH and PCC that for any positive integer $k$ and any fixed $c>0$, the proportion of $k$-consecutive $c$-moderate gaps between ordinates of the zeta zeros is at least $1-k\mc{P}(c)$ (provided this quantity is positive), where $\mc{P}(c)$ is defined in \eqref{Eq:DefinitionPAlpha}. Our next result shows unconditionally that there are positive proportions of $k$-consecutive $c$-moderate gaps for $k=3, 4, 5$ and corresponding values of $c>0$. 

\begin{thm}\label{Thm:Kconcecutive}
Put $c_3=0.1458>\frac17$, $c_4=0.1299>\frac18$ and $c_5=0.1182>\frac19$.  Let  $B_k(T)$ be the number of zeros $\rho_n$ such that $0<\gamma_n\leq T$ and $\min\big(\gamma_{n+k}-\gamma_{n+k-1}, \dots, \gamma_{n+1}-\gamma_n\big)> 2\pi c_k/\log T$. Then we have 
$$\liminf_{T\to \infty} \frac{B_3(T)}{N(T)}\geq 0.345, \qquad \liminf_{T\to \infty} \frac{B_4(T)}{N(T)} \geq 0.1812 \qquad \textup{ and } \qquad \liminf_{T\to \infty} \frac{B_5(T)}{N(T)}\geq 0.0175.
$$
\end{thm}

\begin{rem}\label{Rem:Others}
Independently, with the help of AI,  Ang\'on, Ross, Tang and Zhao \cite{ARTZ} obtained smaller proportions for $k$-consecutive moderate gaps between the zeros of the Riemann zeta function, although the gaps they obtained are larger. More precisely, they showed unconditionally that the proportion of $k$ consecutive $1/4$-moderate gaps between the zeta zeros is at least $34.5\%$ if\footnote{This is equivalent to saying that the proportion of $1/4$-isolated zeros is at least $34.5\%$.} $k=2$, $18.12\%$ if $k=3$, and $1.75\%$ for $k=4$. Their work gives no positive proportion for $k=5$. Their argument combines ideas from the work of the author \cite{La}, together with several spectral ingredients from \cite{AlFu}. 
\end{rem}

To prove our results we further develop the Hilbert space approach via the pair correlation method, which was recently introduced by the author in \cite{La} to give a new conceptually simpler proof of a result discovered by an internal research version of Claude (see \cite{AlFu}) that more than $2/3$ of the zeros of the Riemann zeta function are simple and on the critical line. Our new approach allowed us to discover a ``rigidity'' phenomenon, missed by the Claude proof, which asserts that the properties of a zero being simple or on the critical line are ultimately connected. In particular, we showed in \cite[Theorem 1.1] {La} that the average of the proportions of simple zeros and of zeros on the critical line exceeds $5/6$.  
Here we develop this rigidity phenomenon further, which allows us to show that a large proportion of the zeros on the line, produced by the method, are not only simple but are in fact $c$-isolated for some positive constant $c>0$. To extract such information, our key observation is that zeros on the critical line that are extremely close behave roughly like double zeros, so their contribution to the Hilbert space ``structure'' in Propositions \ref{Pro:IsolatedZeros}, \ref{Pro:KConsecutiveGaps} and \ref{Pro:IsolatedSMALL} below can be ``grouped'' together. This allows us to show that more than half of the zeros are $1/6$-isolated. Our idea can be further pushed  to improve the proportion of these zeros to $>13/20$ when the constant $c$ in Theorem \ref{Thm:Isolated0} is very small, by exploiting the additional rigidity of the extremely close zeros. To prove this last result we also use an explicit zero density estimate due to Conrey \cite{CoBLMS} in order to control the contribution of zeros away from the critical line. Moreover, since the obtained proportion for isolated zeros is very large, a simple combinatorial argument combined with the proof of  \eqref{Thm:Isolated} allows us to investigate $k$-consecutive gaps between the zeros for $k=3,4,5$, leading to Theorem \ref{Thm:Kconcecutive}.

\subsection{Application to zeros of linear combinations of two Dirichlet $L$-functions}\label{Sec:LinearCombinations}
The study of zeros of linear combinations of $L$-functions goes back to the work of Davenport and Heilbronn \cite{DaHe}, who showed in 1936 that certain Epstein zeta functions have infinitely many zeros in the half-plane of absolute convergence. The Epstein zeta function of a quadratic form of class number $\geq 2$ can be expressed as a linear combination of Hecke $L$-functions and as such does not possess an Euler product. This motivated the study of the zeros of general linear combinations of $L$-functions. We refer the reader to \cite{LaLe} for the history of this subject and related results.

Although a linear combination of $L$-functions might have many zeros off the critical line, Montgomery conjectured that a linear combination of primitive $L$-functions satisfying certain natural conditions (notably the same functional equation and a certain notion of independence) has $100\%$ of its non-trivial zeros on the critical line. 
In their beautiful paper \cite{BoHe}, Bombieri and Hejhal  proved Montgomery's conjecture under the assumption that all the involved $L$-functions satisfy the generalized Riemann hypothesis and an analogue of the well spacing hypothesis \eqref{Eq:WellSpacingHyp}. Moreover, in \cite[Theorem C]{BoHe} they also proved conditionally on the GRH alone that a linear combination of two Dirichlet $L$-functions with the same conductor and parity has at least $34.5\%$ of its non-trivial zeros on the critical line.  Unconditionally, Selberg \cite{Sel98} proved in 1998 that a positive proportion of the zeros of any linear combination of $J\geq 2$ Dirichlet $L$-functions attached to primitive characters of the same conductor and parity lie on the critical line. In his thesis, Dousselin \cite{Do} refined Selberg's method and made it explicit. In particular, he showed that in the special case $J=2$, this proportion is at least $7.38\times 10^{-9}$. We refine the Bombieri--Hejhal method, making it unconditional, and combine it with our results on isolated zeros (which we extend to Dirichlet $L$-functions) to obtain a substantial improvement of this proportion. 
\begin{thm}\label{Thm:ZerosLinearComb}
Let $b_1$, $b_2$ be non-zero real numbers. Let $q\geq 3$ be a positive integer, and $\chi_1$, $\chi_2$ be distinct primitive characters modulo $q$, such that $\chi_1(-1)=\chi_2(-1)$. For $j=1,2$ let $\theta_j$ be such that\footnote{The phase $\theta_j$ is chosen in order to make the function $e^{i\theta_j}\Lambda(s, \chi_j)$ real on the critical line, where $\Lambda(s, \chi_j)$ is the completed $L$-function attached to $\chi_j$, see the functional equation \eqref{Eq:FunctionEquation}.} $e^{2i\theta_j}=i^{a}\overline{\tau(\chi_j)}/\sqrt{q}$, where $a=(1-\chi_j(-1))/2$ and $\tau(\chi_j)$ is the Gauss sum attached to $\chi_j$. Then at least $31\%$ of the non-trivial zeros of the linear combination $b_1 e^{i\theta_1} L(s, \chi_1)+ b_2 e^{i\theta_2} L(s, \chi_2)$ lie on the critical line and have odd multiplicity. 
\end{thm}
\begin{rem}\label{Rem:ZerosLinearComb}
Our unconditional refinement of the Bombieri--Hejhal method (see Proposition \ref{Pro:BombieriHejhal} below) works only if one can obtain proportions $\geq\delta_j$, for $j=1, 2$, for $c$-isolated zeros (with $c>0$ small and fixed) of the Dirichlet $L$-functions $L(s, \chi_j)$, such that $\delta_1+\delta_2>1.$ In this case we can show that the proportion of zeros on the critical line of the linear combination $b_1 e^{i\theta_1} L(s, \chi_1)+ b_2 e^{i\theta_2} L(s, \chi_2)$ is $\geq  \delta_1+\delta_2-1$.  In particular since we prove the analogue of \eqref{Thm:IsolatedSmall}, showing that at least $0.6552$ of the zeros of a Dirichlet $L$-function are $c$-isolated for some $c>0$, we obtain the proportion $\delta_0=2\left(0.6552-1/2\right)=0.3104$ for the zeros of the linear combination $b_1 e^{i\theta_1} L(s, \chi_1)+ b_2 e^{i\theta_2} L(s, \chi_2)$ on the critical line. 
\end{rem}
\subsection*{Acknowledgments} The author is supported by a junior chair of the Institut Universitaire de France. Parts of this work (namely the estimate \eqref{Thm:Isolated} and Theorem \ref{Thm:Kconcecutive}) were completed while the author was visiting the Max Planck Institute for Mathematics in Bonn (during September 2026), which he thanks for its warm hospitality and excellent working conditions. 

This paper contains no AI generated text or ideas.


\section{Hilbert space inequalities for the number of consecutive gaps in a general multiset of complex numbers}

\subsection{Isolated real elements of a multiset}
For a compactly supported function $f\in L^{1}(\mathbb R)$ we define its Fourier transform $\widehat{f}:\C\to \C$ by
\begin{equation}\label{hat-g}
\widehat{f}(\xi) = \int_{-\infty}^{\infty} f(u)e^{-2\pi i\xi u }\,du, \quad \text{ for all } \xi\in \C.
\end{equation}
Note that this is an entire function since $f$ is compactly supported. 

Let $\lambda>0$ be a real number and $\eta\in L^2(\R)$ be a real valued even function with $\textup{supp} (\eta) \subset (-\lambda, \lambda) $, such that $\widehat{\eta^2}(0)=1$. We define the kernel
$$K(\xi):= \widehat{\eta^2}(\xi), \quad \text{ for } \quad \xi\in \C.$$
Since $\eta$ is even and real valued we have $K(\xi)\in \R$ for all $\xi\in \R$.

Throughout this section we let $\mc{Z}$ be a non-empty finite multiset of complex numbers which is \emph{invariant under complex conjugation}, that is for all $z\in \mc{Z}$ we have $\overline{z}\in \mc{Z}$ and $m_{\overline{z}}=m_z,$ where $m_z$ denotes the multiplicity of $z$ in $\mc{Z}$. Unless otherwise specified, sums over a multiset are understood to count elements with their multiplicities. We prove the following proposition.
\begin{pro}\label{Pro:IsolatedZeros}
Let $\tau>0$ be a real number and put $\kappa_0=\min \{K(x): x\in [-2\tau, 2\tau]\}$. Assume that $\kappa_0\geq 5/6$. 
Then we have 
$$
\Big|\big\{x\in \mc{Z}\cap \R: |x-\re(z)|>\tau \text{ for all } z\in \mc{Z}\setminus\{x\}\big\}\Big|\geq \frac52\sum_{z\in \mc{Z}} 1 -
 \frac32\sum_{z,s\in\mathcal Z}K(z-s)^2. 
$$
\end{pro}
Note that for $x\in  \mc{Z}\cap \R$,  the condition $|x-\re(z)|>\tau \text{ for all } z\in \mc{Z}\setminus\{x\}$ implies that $x$ is simple. 

\begin{proof}[Proof of Proposition \ref{Pro:IsolatedZeros}]

We start the proof similarly to that of Proposition 2.1 of \cite{La}. For all $z\in \C$ we define the functions $f_z, g_z, h_z:\R\to \C$ by
\begin{equation}\label{Eq:DefinitionFz}
 f_z(u):=\eta(u) e^{-2\pi i u z}, \quad 
 g_{z}(u) := \frac{f_z(u)+f_{\overline{z}}(u)}{2}, \quad \text{ and } \quad h_{z}(u)= \frac{f_z(u)-f_{\overline{z}}(u)}{2i}.
\end{equation}
We record the following facts established in the proof of Proposition 2.1 of \cite{La}. First we have
\begin{equation}\label{Eq:KernelFactorization}
K(z-\overline{s})=\int_{-\lambda}^{\lambda} f_z(u)\overline{f_s(u)}\,du,
\end{equation}
from which we deduce that for all $x\in \R$ and $z\in \C$ we have
\begin{equation}\label{Eq:IntegralFz}
\int_{-\lambda}^{\lambda} |f_x(u)|^2du=1, \quad \textup{ and } \quad \int_{-\lambda}^{\lambda} |g_z(u)|^2 du- \int_{-\lambda}^{\lambda} |h_z(u)|^2 du =1.
\end{equation}
Now for all $(u, v) \in \R^2$ we define
$$
F(u,v):=\sum_{z\in\mathcal Z}f_z(u)f_z(v).
$$
Then we have
\begin{equation}\label{Eq:SecondMoment}
\sum_{z, s\in \mc{Z}} K^2 (z-s)
=\int_{-\lambda}^{\lambda} \int_{-\lambda}^{\lambda} |F(u, v)|^2 \,du \,dv.
\end{equation}
Let $\mc{R}_1$ be the set of simple real elements of $\mc{Z}$, $\mc{R}_2$ be the set of real elements of $\mc{Z}$ with multiplicity $\geq 2$, and $\mc{S}$ be the set of its distinct non-real elements.  Since $\mc{Z}$ is invariant under complex conjugation we can list the elements of these sets as 
\begin{equation}\label{Eq:ListElementsZ}
\mc{R}_1= \{x_1, \dots, x_n\}, \quad \mc{R}_2= \{x_{n+1}, \dots, x_{n+r}\}, \ \text{ and } \ \mc{S}=\{z_1, \overline{z}_1, \dots, z_q, \overline{z}_q\}.
\end{equation}
Then we can write 
\begin{equation}\label{Eq:IdentityFUV}
F(u, v) 
= \sum_{\ell=1}^{n+ r} m_{x_\ell}f_{x_\ell}(u)f_{x_\ell}(v)+ 2\sum_{\ell=1}^q  m_{z_\ell} \left(g_{z_\ell}(u) g_{z_\ell}(v)- h_{z_\ell}(u) h_{z_\ell}(v)\right).
\end{equation}

We now deviate from the proof of \cite[Proposition 2.1]{La} and consider a different construction of the nested vector spaces. To this end we split the set of simple real elements $\mc{R}_1$ into two parts: the ``isolated elements'' and the ``clustered elements''. Here the set of isolated points is defined by 
\begin{equation}\label{Eq:DefinitionIsolated}
\mc{I}= \big\{x\in \mc{R}_1 \textup{ such that } |x-y|>2\tau \textup{ for all } y\in \mc{R}_1\setminus\{x\}\big\}.
\end{equation}
These are the elements of $\mc{R}_1$ that are at least $2\tau$ apart from all other elements of $\mc{R}_1.$ 
Let $d=|\mc{I}|$. We can rearrange the elements of $\mc{R}_1$ so that 
$$ \mc{I}= \big\{x_1, \dots, x_d\}. $$
Note that for all $z\in \mc{R}_2\cup \mc{S}$, $|\re(z)-x|\leq \tau$ for at most one element $x\in \mc{I}$. Therefore, the quantity we seek to bound is 
\begin{equation}\label{Eq:QuantityConsecutive}
\Big|\big\{x\in \mc{Z}\cap \R: |x-\re(z)|>\tau \text{ for all } z\in \mc{Z}\setminus\{x\}\big\}\Big|\geq |\mc{I}|- |\mc{R}_2|-\frac{|\mc{S}|}{2}= d-(r+q).
\end{equation}

Let $\mc{C}=\mc{R}_1\setminus \mc{I}$ be the clustered elements. We can split $\mc{C}$ into a disjoint union of maximal clusters of consecutive elements which are $\leq 2\tau$ apart, namely
\begin{equation}\label{Eq:DefinitionClusters}
\mc{C}= \mc{C}_1 \cup \mc{C}_2 \cup \cdots \cup \mc{C}_t,
\end{equation}
where each cluster has the form
$\mc{C}_j=\{y_1, \dots, y_e\},$
with $0\leq y_{\ell+1}-y_{\ell}\leq 2\tau$. We group the elements of each cluster into pairs of consecutive elements if the cardinality of the cluster is even, and then into pairs and one remaining triple if its cardinality is odd. 

Our strategy consists in grouping the contributions in the Hilbert space structure of the functions $f_{x}$ for pairs and triples that are nearby. The heuristic reason behind this is that two distinct elements $y_1$ and $y_2$ that are extremely close behave roughly like an element with multiplicity $2$. This suggests to replace\footnote{Note that if $f_{y_2}=f_{y_1}$, then $f_{y_1}+f_{y_2}=2f_{y_1}$ which is a generator of the vector space $\textup{Span}(f_{y_1})$.} the vectors $f_{y_1}$ and $f_{y_2}$ by the vector $f_{y_1}+f_{y_2}.$
Let $\mc{G}_1, \dots, \mc{G}_a$ be all the pairs that are obtained from all such clusters, and let $\mc{T}_1, \dots, \mc{T}_b$ be all the triples.  
 Then we have 
\begin{equation}\label{Eq:IdentitySimpleReal}
n= d+ 2a+ 3b. 
\end{equation}
 For each pair $\mc{G}_{\ell}$ and each triple $\mc{T}_{j}$ we define
\begin{equation}\label{DefPairsTriples}
 P_{\ell}:= P_{\mc{G}_\ell}= \sum_{y\in \mc{G}_{\ell}}f_{y}, \quad \textup{ and } \quad
 Q_{j}:= Q_{\mc{T}_j}= \sum_{y\in \mc{T}_{j}}f_{y}.
\end{equation} 
We consider the Hilbert space $L^2((-\lambda, \lambda), \C)$ endowed with its usual inner product 
$$\langle f, g\rangle = \int_{-\lambda}^{\lambda} f(u)\overline{g(u)}du,$$
and define the following two nested finite dimensional subspaces of $L^2((-\lambda, \lambda), \C)$:
$$
U=\textup{Span}\big(P_1, \dots, P_a, Q_1, \dots, Q_b, f_{x_{n+1}}, \dots, f_{x_{n+r}}, g_{z_1}, \dots, g_{z_q}\big), 
$$ and 
$$ V=\textup{Span}\big(f_{x_1}, \dots, f_{x_d}, P_1, \dots, P_a, Q_1, \dots, Q_b, f_{x_{n+1}}, \dots, f_{x_{n+r}}, g_{z_1}, \dots, g_{z_q}\big).$$
We also put 
$$D_U:=\dim U, \quad  \text{ and } \ D_V:= \dim V. $$
 Note that for distinct complex numbers $z_1, \dots, z_{\ell}$, the functions $e^{z_1 u}, e^{z_2 u}, \dots, e^{z_\ell u}$ are linearly independent over $\C$. Therefore, since the set on which $\eta$ is non-zero has positive measure (because $\widehat{\eta^2}(0)=1$), it follows that $f_{x_1}, \ldots, f_{x_{n+r}}, g_{z_1}, \ldots, g_{z_q}$ are linearly independent over $\C$. Therefore, $f_{x_1}, \dots, f_{x_d}, P_1, \dots, P_a, Q_1, \dots, Q_b, f_{x_{n+1}}, \dots, f_{x_{n+r}}, g_{z_1}, \dots, g_{z_q}$ are also linearly independent over $\C$ and we get 
\begin{equation}\label{Eq:DimensionCalculations}
  D_U= a+b+r+q \quad \text{ and } \ D_V= d+a+b+ r+q.  
 \end{equation}
Note that all of the functions $f_{x_{\ell}}, g_{z_{\ell}}$ satisfy the property $
\overline{\phi(u)}= \phi(-u)$ for all $u\in \mathbb{R}$, and hence the functions $P_\ell$, $Q_\ell$ satisfy this property as well. Therefore, for any  pair of functions $\phi_1, \phi_2\in 
\{f_{x_1}, \dots, f_{x_d}, P_1, \dots, P_a, Q_1, \dots, Q_b, f_{x_{n+1}}, \dots, f_{x_{n+r}}, g_{z_1}, \dots, g_{z_q}\}$
we have $\overline{\langle \phi_1, \phi_2\rangle} 
=\langle \phi_1, \phi_2\rangle$ and hence all inner products between these functions are real.
We now apply the standard Gram-Schmidt process to construct an orthonormal basis $(\psi_1, \dots, \psi_{D_V})$ of $V$, such that  $(\psi_1, \dots, \psi_{D_U})$ is a basis for $U$. By the above discussion, all the Gram–Schmidt coefficients are real and hence $\overline{\psi_j(u)}=\psi_j(-u)$ for all $u\in \mathbb{R}$. Therefore, all of the inner products appearing below are real.

Since the functions $\psi_j(u)\psi_j(v)$ for $1\leq j\leq D_V$ are orthonormal in $L^2\big((-\lambda, \lambda)^2,\C\big)$, it follows from Bessel's inequality that
\begin{equation}\label{Eq:OrthogonalInequality}
\int_{-\lambda}^{\lambda}\int_{-\lambda}^{\lambda} |F(u, v)|^2 du dv \geq \sum_{j=1}^{D_V}|\alpha_j|^2,
\end{equation}
where 
\begin{align}\label{Eq:DefinitionAlpha}
\alpha_j&:= \int_{-\lambda}^{\lambda}\int_{-\lambda}^{\lambda} F(u, v) \overline{\psi_j(u)\psi_j(v)}\,du \,dv \nonumber\\
&= \sum_{\ell=1}^{n} \langle f_{x_{\ell}}, \psi_j\rangle ^2+ \sum_{\ell=n+1}^{n+ r} m_{x_{\ell}}\langle f_{x_{\ell}}, \psi_j\rangle ^2 + 2\sum_{\ell=1}^q  m_{z_{\ell}} \left(\langle g_{z_{\ell}}, \psi_j\rangle ^2- \langle h_{z_{\ell}}, \psi_j\rangle ^2\right) 
\end{align}
by \eqref{Eq:IdentityFUV}. We use the basic inequalities $y^2\geq 2y-1$ and $y^2\geq 4y-4$, valid for all $y\in \R$ to obtain 
\begin{equation}\label{Eq:SqaureALPHA_j}
\sum_{j=1}^{D_U} \alpha_j^2\geq 4 \sum_{j=1}^{D_U} \alpha_j - 4D_U, \quad \textup{ and } \sum_{j=D_U+1}^{D_V} \alpha_j^2\geq 2 \sum_{j=D_U+1}^{D_V} \alpha_j - (D_V-D_U).
\end{equation}
Summing these inequalities and combining them with \eqref{Eq:SecondMoment} and \eqref{Eq:OrthogonalInequality} implies that 
\begin{equation}\label{Eq:MainInequalityKernel}
 \sum_{z, s\in \mc{Z}} K^2 (z-s)\geq 2 \sum_{j=1}^{D_U} \alpha_j+  2\sum_{j=1}^{D_V} \alpha_j -3D_U-D_V.
\end{equation}
To finish the proof, we need to obtain a lower bound on the sums $\sum_{j=1}^{D_U} \alpha_j$ and $\sum_{j=1}^{D_V} \alpha_j$. First, by Bessel's inequality we have 
\begin{equation}\label{Eq:Bessel}
\sum_{j=1}^{D_V}  \langle h_{z_{\ell}}, \psi_j\rangle ^2 \leq \int_{-\lambda}^{\lambda} |h_{z_{\ell}}(u)|^2 du.
\end{equation}
Moreover, since $f_{x_{n+1}}, \dots, f_{x_{n+r}}, g_{z_1}, \dots, g_{z_{q}} \in U\subset V$, it follows from Parseval's theorem, \eqref{Eq:IntegralFz}  and \eqref{Eq:Bessel} that
\begin{align}\label{Eq:EasyPartParsevalBessel}
    &\sum_{j=1}^{D_U} \left(\sum_{\ell=n+1}^{n+ r} m_{x_{\ell}}\langle f_{x_{\ell}}, \psi_j\rangle ^2 + 2\sum_{\ell=1}^q  m_{z_{\ell}} \left(\langle g_{z_{\ell}}, \psi_j\rangle ^2- \langle h_{z_{\ell}}, \psi_j\rangle ^2\right)\right) \nonumber\\
    &\geq  \sum_{\ell=n+1}^{n+ r} m_{x_{\ell}} \int_{-\lambda}^{\lambda} |f_{x_{\ell}}(u)|^2 du+ 2\sum_{\ell=1}^q  m_{z_{\ell}} \left(\int_{-\lambda}^{\lambda} |g_{z_{\ell}}(u)|^2 du-\int_{-\lambda}^{\lambda} |h_{z_{\ell}}(u)|^2 du \right)\\
    &= \sum_{\ell=n+1}^{n+ r} m_{x_{\ell}} + 2\sum_{\ell=1}^q  m_{z_{\ell}}.\nonumber
\end{align}
Combining this estimate with \eqref{Eq:DefinitionAlpha} gives with $W$ being either $U$ or $V$
\begin{equation}\label{Eq:LowerBoundSUMAlphaBoth}
    \sum_{j=1}^{D_W} \alpha_j \geq \sum_{\ell=1}^{n}\sum_{j=1}^{D_W} \langle f_{x_{\ell}}, \psi_j\rangle ^2 + \sum_{\ell=n+1}^{n+ r} m_{x_{\ell}} + 2\sum_{\ell=1}^q  m_{z_{\ell}}. 
\end{equation}
We first handle the contribution of the isolated elements $\{x_1, \dots, x_d\}$. Since $f_{x_1}, \dots, f_{x_d}\in V$ we have 
\begin{equation}\label{Eq:PaservalLastIsolated}
  S_0:=\sum_{\ell=1}^{d}\sum_{j=1}^{D_U} \langle f_{x_{\ell}}, \psi_j\rangle^2\geq 0 \quad \text{ and } \quad    \sum_{\ell=1}^{d}\sum_{j=1}^{D_V} \langle f_{x_{\ell}}, \psi_j\rangle ^2= \sum_{\ell=1}^{d} \int_{-\lambda}^{\lambda} |f_{x_\ell}(u)|^2du= d,
\end{equation}
where the second identity follows from Parseval's theorem and \eqref{Eq:IntegralFz}. 

It remains to handle the contribution of the clusters. Here we recall that $\{x_{d+1}, \dots, x_n\}$ is the disjoint union of the pairs $\{\mc{G}_\ell\}_{1\leq \ell\leq a}$ and the triples $\{\mc{T}_\ell\}_{1\leq \ell\leq b}$. For each pair $\mc{G}_\ell$ we have by the Cauchy-Schwarz inequality that 
$$
\sum_{y\in \mc{G}_{\ell}} \langle f_{y}, \psi_j\rangle^2 \geq \frac{1}{2} \bigg(\sum_{y\in \mc{G}_{\ell}} \langle f_{y}, \psi_j\rangle\bigg)^2= \frac12\Big\langle \sum_{y\in \mc{G}_\ell}f_{y}, \psi_j\Big\rangle^2= \frac12\langle P_\ell, \psi_j\rangle^2.
$$
Similarly one has 
$$
\sum_{y\in \mc{T}_{\ell}} \langle f_{y}, \psi_j\rangle^2 \geq \frac{1}{3} \bigg(\sum_{y\in \mc{T}_{\ell}} \langle f_{y}, \psi_j\rangle\bigg)^2= \frac13\langle Q_\ell, \psi_j\rangle^2.
$$
Hence, since $P_1, \dots, P_a, Q_1, \dots, Q_b\in U \subset V$ we deduce from Parseval's theorem that 
\begin{align}\label{Eq:ParsevalPQ1}
\sum_{\ell=d+1}^n \sum_{j=1}^{D_V} \langle f_{x_{\ell}}, \psi_j\rangle ^2\geq \sum_{\ell=d+1}^n \sum_{j=1}^{D_U} \langle f_{x_{\ell}}, \psi_j\rangle ^2 &= \sum_{\ell=1}^a \sum_{j=1}^{D_U} \sum_{y\in \mc{G}_{\ell}} \langle f_{y}, \psi_j\rangle ^2+ \sum_{\ell=1}^b \sum_{j=1}^{D_U} \sum_{y\in \mc{T}_{\ell}} \langle f_{y}, \psi_j\rangle ^2\nonumber\\
& \geq \frac12\sum_{\ell=1}^a \sum_{j=1}^{D_U} \langle P_\ell, \psi_j\rangle^2 +   \frac13\sum_{\ell=1}^b \sum_{j=1}^{D_U} \langle Q_\ell, \psi_j\rangle^2 \\
& = \frac12 \sum_{\ell=1}^a \int_{-\lambda}^{\lambda} |P_\ell(u)|^2 du  +   \frac13\sum_{\ell=1}^b \int_{-\lambda}^{\lambda} |Q_\ell(u)|^2 du.\nonumber  
\end{align}
We now derive lower bounds for the $L_2$-norms of $P_\ell$ and $Q_\ell$, which will replace \eqref{Eq:IntegralFz}. By definition, we can write $P_{\ell}= f_{y_1}+ f_{y_2}$ for some real elements $y_1\leq y_2$ such that $y_2-y_1\leq 2\tau$. Therefore, by \eqref{Eq:KernelFactorization} and \eqref{Eq:IntegralFz} we have 
\begin{equation}\label{Eq:L2NormP}
\int_{-\lambda}^{\lambda}|P_{\ell}(u)|^2 du= 2+ 2 \int_{-\lambda}^{\lambda}f_{y_1}(u)\overline{f_{y_2}(u)} du=  2+ 2 K(y_2-y_1)\geq 2+2\kappa_0,
\end{equation}
by the definition of $\kappa_0$. 
Similarly we can write $Q_{\ell}= f_{w_1}+ f_{w_2}+f_{w_3}$ for some real elements $w_1\leq w_2\leq w_3$ such that $w_2-w_1\leq 2\tau$ and $w_3-w_2\leq 2\tau$. Expanding the square and using \eqref{Eq:KernelFactorization} and \eqref{Eq:IntegralFz} we get 
\begin{equation}\label{Eq:L2NormQ}
\int_{-\lambda}^{\lambda}|Q_{\ell}(u)|^2 du=  3+ 2 K(w_2-w_1)+ 2K(w_3-w_2)+ 2K(w_3-w_1)\geq 3+4\kappa_0+ 2\kappa_1,
\end{equation} 
where $\kappa_1:=\min \{K(x): x\in [-4\tau, 4\tau]\}$. 

We now combine the estimates \eqref{Eq:LowerBoundSUMAlphaBoth}--\eqref{Eq:L2NormQ} to obtain
\begin{equation}\label{Eq:LowerAlphaU}
  \sum_{j=1}^{D_U} \alpha_j  \geq S_0+ a(1+\kappa_0)+ b\left(1+\frac43 \kappa_0+ \frac23 \kappa_1\right)+ \sum_{\ell=n+1}^{n+ r} m_{x_{\ell}} + 2\sum_{\ell=1}^q  m_{z_{\ell}},
\end{equation}
and 
\begin{equation}\label{Eq:LowerAlphaV}
    \sum_{j=1}^{D_V} \alpha_j \geq d+ a(1+\kappa_0)+ b\left(1+\frac43 \kappa_0+ \frac23 \kappa_1\right)+  \sum_{\ell=n+1}^{n+ r} m_{x_{\ell}} + 2\sum_{\ell=1}^q  m_{z_{\ell}}. 
\end{equation}
To simplify our notation we let $N=\sum_{z\in \mc{Z}} 1.$ Then we observe that 
\begin{equation}\label{Eq:EstimateN}
N= d+ 2a+3 b+ \sum_{\ell=n+1}^{n+r} m_{x_{\ell}} + 2\sum_{\ell=1}^q  m_{z_{\ell}}\geq d+2a+3b+2(r+q).
\end{equation}
Combining the first identity with \eqref{Eq:DimensionCalculations}, \eqref{Eq:MainInequalityKernel}, \eqref{Eq:LowerAlphaU} and \eqref{Eq:LowerAlphaV} we get
\begin{equation}\label{Eq:MainInequalityKernel2} \sum_{z, s\in \mc{Z}} K^2 (z-s)\geq 4N-3d-4r-4q-4a(2-\kappa_0)-4b\Big(3-\frac43 \kappa_0-\frac23 \kappa_1\Big) +2S_0.
\end{equation}
By the Cauchy-Schwarz inequality and the fact that $\eta$ is even we have for all $u\in \R$
\begin{align}\label{Eq:CauchySchwarzKernel} K(u)^2 = 4 \left(\int_0^{\lambda} \eta^2(v)\cos(2\pi u v) dv\right)^2&\leq 4 \left(\int_{0}^\lambda \eta^2(v)dv\right)\left(\int_{0}^\lambda \eta^2(v) \cos(2\pi u v)^2 dv\right)\nonumber\\
& = \int_{0}^\lambda \eta^2(v)\big(1+\cos(4\pi u v)\big) dv= \frac{1+K(2u)}{2},
\end{align}
since $2\int_0^{\lambda}\eta(u)^2du=1$. Thus, using our assumption $\kappa_0\geq 5/6$ we obtain $\kappa_1\geq 2\kappa_0^2-1\geq 7/18$.

By \eqref{Eq:QuantityConsecutive} it suffices to prove a lower bound for $d-(r+q).$ To use the inequality on the right hand side of \eqref{Eq:EstimateN} we need to express the term $ 3d+4r+4q$ that appears on the right hand side of \eqref{Eq:MainInequalityKernel2} as a linear combination of $d-(r+q)$ and $d+2(r+q)$. An easy calculation gives $3d+4r+4q= \frac23(d-(r+q))+ \frac7{3}(d+2(r+q)).$ Therefore, it follows from \eqref{Eq:EstimateN} and the fact that $S_0\geq 0$ that the sum on the left hand side of \eqref{Eq:MainInequalityKernel2} is
\begin{align*}
&\geq 4N-\frac23(d-(r+q))-\frac7{3} (d+2r+2q)- 4a(2-\kappa_0)- 4b\Big(3-\frac43 \kappa_0-\frac23 \kappa_1\Big)\\
&\geq 4N-\frac23(d -(r+q))- \frac73\big(d+2a+3b+2(r+q)\big)\geq \frac{5}{3} N- \frac23(d -(r+q)). 
\end{align*}
Hence we deduce that 
$$ d- (r+q) \geq \frac52 N- \frac32  \sum_{z, s\in \mc{Z}} K^2 (z-s).$$
Inserting this estimate in \eqref{Eq:QuantityConsecutive} completes the proof. 

\end{proof}

\subsection{The case of $k$-consecutive gaps}
Let $k\geq 2$ be an integer and define a $k$-block to be a block of $k$ gaps between $k+1$ consecutive points. We say that a $k$-block has width $> \delta$ if all gaps between the elements forming this block are $> \delta$.  Our goal is to count the number of $k$-blocks of width $> \tau$ formed by the real parts of the elements of $\mc{Z}$. We call such a block ``good''. Note that a good block will have all its interior points necessarily simple and real since only the two end points might be non-simple or non-real. We will first ignore at most $O_k(1)$ starting and ending blocks, in order to assume that every point appears in exactly $k+1$ of the $k$-blocks. To simplify our notation we define $\re(\mc{Z})$ to be the multiset of real parts of elements of $\mc{Z}$. We prove the following proposition.

\begin{pro}\label{Pro:KConsecutiveGaps}
Let $k\geq 3$ be a positive integer. Let $\tau>0$ be a real number and put $\kappa_0=\min \{K(x): x\in [-2\tau, 2\tau]\}$ and $\kappa_1=\min \{K(x): x\in [-4\tau, 4\tau]\}$. Assume that 
$
\kappa_0\geq 1-1/(2(k+1))$. 
Then the number of $k$-blocks of width $>\tau$ formed by the elements of $\re(\mc{Z})$ is 
$$ \geq \frac{k+3}{2}\sum_{z\in \mc{Z}}1- \frac{k+1}{2}\sum_{z, s\in \mc{Z}} K^2 (z-s)+O_k(1). $$
\end{pro}
\begin{proof}
As in the proof of Proposition \ref{Pro:IsolatedZeros} we let $\mc{I}$ be the set of isolated real simple points defined in \eqref{Eq:DefinitionIsolated}, and consider the maximal clusters defined in \eqref{Eq:DefinitionClusters} whose number is $t$. We now collapse each of these maximal clusters into a single ``point'' and call these ``points'' $c_1, \dots, c_t$.
Thus the sequence $\mc{K}=\{x_1, \dots, x_d, c_1, \dots, c_t\}$ is $2\tau$ well-spaced (meaning that the distance between every two elements is $> 2\tau$). Now we insert the real parts of the elements in $\mc{R}_2\cup \mc{S}$. Here, it is understood that we collapse a real point of multiplicity $m$ in $\mc{Z}$ into one point in $\mc{R}_2$, and we do the same for the real parts of the non-real complex conjugate elements of $\mc{Z}$. Moreover, if the inserted point lands inside a cluster it becomes part of it. Then one can see that any newly inserted element $y$ can be within $\tau$ of at most one element $x\in \mc{K}$. If this is the case we collapse the pair $\{x, y\}$ into a new cluster.  
Let $L$ be the number of newly formed collapsed clusters which we call $e_1, \dots, e_L$. Then we have $L\leq t+r+q$.

Let $\mc{I}_1$ be the set of those $x\in \mc{I}$ such that $|x-\re(z)|>\tau$ for all $z\in \mc{Z}\setminus\{x\}$. We write $\mc{I}_1= \{w_1, \dots, w_J\}$. By \eqref{Eq:QuantityConsecutive} we have $J\geq d-r-q$. We now consider the sequence $\{w_1, \dots, w_J, e_1, \dots e_L\}$. The number of $k$-blocks formed by elements of this sequence is $J+L+O_k(1)$.    
Furthermore, with our convention that every cluster $e_j$ is collapsed into a single ``point'', such a point can be the interior point of at most $k-1$ of the $k$-blocks of the elements $w_1, \dots, w_J, e_1, \dots, e_L$. We will exclude these $k$-blocks, so that the remaining ones have all their interior elements among the $w_j$'s and hence are good $k$-blocks. Here, it is understood that a block starting at a collapsed cluster starts at its rightmost original point, and a block ending at a collapsed cluster ends at its leftmost original point. Therefore, we deduce that the number of good $k$-blocks formed by the elements of $\re(\mc{Z})$ is 
\begin{align}\label{Eq:FinalCountConsecutiveGaps}
&\geq J+L- (k-1)L+ O_k(1)= J-(k-2)L+O_k(1) \nonumber\\
&\geq d-(k-1)(r+q)- (k-2)t+ O_k(1). 
\end{align}

Now, in order to obtain a lower bound for this quantity, we will proceed as before to first express the term $3d+4r+4q$ on the right hand side of \eqref{Eq:MainInequalityKernel2} as a linear combination of the terms $d-(k-1)(r+q)$ (which appears in \eqref{Eq:FinalCountConsecutiveGaps}) and $d+2(r+q)$ (which appears on the right hand side of \eqref{Eq:EstimateN}). An easy calculation gives 
$$
3d+ 4r+ 4q= \frac{2}{k+1}\big(d-(k-1)(r+q)\big)+ \left(3-\frac{2}{k+1}\right)\big(d+2(r+q)\big).
$$ Note that $t\leq  a+b$. Therefore, we deduce from \eqref{Eq:EstimateN} and \eqref{Eq:MainInequalityKernel2} that 
\begin{align*}
&\sum_{z, s\in \mc{Z}} K^2 (z-s)\geq 4N-\frac{2}{k+1}\big(d-(k-1)(r+q)-(k-2)t\big)- \left(3-\frac{2}{k+1}\right)\big(d+2(r+q)\big)\\
& \qquad \qquad  - a\left(8-4\kappa_0+\frac{2k-4}{k+1}\right)-b\Big(12-\frac{16}{3} \kappa_0-\frac83 \kappa_1+ \frac{2k-4}{k+1}\Big)\\
&\geq  4N-\frac{2}{k+1}\big(d-(k-1)(r+q)-(k-2)t\big)- \left(3-\frac{2}{k+1}\right)\big(d+2(r+q)+2a+3b\big)\\
& \geq \frac{k+3}{k+1} N- \frac{2}{k+1}\big(d-(k-1)(r+q)-(k-2)t\big),
\end{align*}
where we have used that $\kappa_0\geq 1-1/(2(k+1))$ by our assumption and $2\kappa_0+\kappa_1\geq 15/8$ by \eqref{Eq:CauchySchwarzKernel}.
This implies 
$$d-(k-1)(r+q)-(k-2)t \geq \frac{k+3}{2}N- \frac{k+1}{2}\sum_{z, s\in \mc{Z}} K^2 (z-s).$$
Combining this estimate with \eqref{Eq:FinalCountConsecutiveGaps} completes the proof.  
\end{proof}

\subsection{The case where $\tau$ is very small} In this section we show that one can improve Proposition \ref{Pro:IsolatedZeros} if $\tau$ is very small.  Indeed we prove
\begin{pro}\label{Pro:IsolatedSMALL}
Let $\tau>0$ be a small real number. Then we have 
\begin{align*}
& \Big|\big\{x\in \mc{Z}\cap \R: |x-\re(z)|>\tau \text{ for all } z\in \mc{Z}\setminus\{x\}\big\}\Big|\\
& \geq 2N -\sum_{z, s\in \mc{Z}} K^2 (z-s) -\sum_{\ell=1}^q \max\left(1-\frac{4K\big(i\im(z_{\ell})\big)^2}{1+K\big(2 i\im(z_{\ell})\big)}, 0\right) +O\left(\tau^2 N\right).  
\end{align*}
\end{pro}
\begin{proof}
First note that for all $x$ with $|x|\leq 4\tau$ we have
\begin{equation}\label{Eq:ApproxKernel0}
 K(x)=\int_{-\lambda}^\lambda \eta^2(u) e^{-2\pi i ux} du= \int_{-\lambda}^\lambda \eta^2(u) (1-2\pi i ux+ O(u^2x^2)\big)du= 1+O(\tau^2),
\end{equation}
since $K(0)=1$ and $\eta$ is even. Thus in this case both $\kappa_0$ and $\kappa_1$ are equal to $1+O(\tau^2)$. Therefore, with the same notation as in the proof of Proposition \ref{Pro:IsolatedZeros} we obtain from \eqref{Eq:EstimateN} and \eqref{Eq:MainInequalityKernel2} that 
\begin{align}\label{Eq:MainInequalitySmallTau}
\sum_{z, s\in \mc{Z}} K^2 (z-s)
&\geq 4N-3d-4r-4q-(4+O(\tau^2))(a+b) +2S_0 \nonumber\\
&\geq  2N-d+2S_0 +O(\tau^2 N).
\end{align}
We will now estimate the contribution of the extra term 
$$ S_0=\sum_{\ell=1}^{d}\sum_{j=1}^{D_U} \langle f_{x_{\ell}}, \psi_j\rangle^2,$$
which we discarded in the proof of Proposition \ref{Pro:IsolatedZeros}. This will be used to control the contribution of the real elements $x_\ell\in \mc{I}$ (where $\mc{I}$ is defined in \eqref{Eq:DefinitionIsolated}) which are extremely close to other elements in $\mc{R}_2$ or to the real part of an element in $\mc{S}$.  Let $1\leq \ell \leq d$, and assume first that there exists $x_m$ with $n+1\leq m\leq n+r$  such that $|x_{\ell}-x_m|\leq \tau$. Write $f_{x_m}\in U$ in the basis $(\psi_1, \dots, \psi_{D_U})$ as 
$$ f_{x_m}= \sum_{j=1}^{D_U} \langle f_{x_m}, \psi_j\rangle \psi_j. 
$$
Since $|x_\ell-x_m|\leq \tau$, it follows from \eqref{Eq:ApproxKernel0} that 
$$ 1+O(\tau^2)= K(x_\ell-x_m)= \langle f_{x_\ell}, f_{x_m}\rangle
=\sum_{j=1}^{D_U}\langle f_{x_m}, \psi_j\rangle\langle f_{x_\ell}, \psi_j\rangle. $$
On the other hand, by the Cauchy-Schwarz inequality we get
\begin{equation}\label{Eq:CSSmallDistance}
\bigg|\sum_{j=1}^{D_U}\langle f_{x_m}, \psi_j\rangle\langle f_{x_\ell}, \psi_j\rangle\bigg|^2 \leq \left(\sum_{j=1}^{D_U}\langle f_{x_\ell}, \psi_j\rangle^2\right) \left(\sum_{j=1}^{D_U}\langle f_{x_m}, \psi_j\rangle^2\right)=  \sum_{j=1}^{D_U}\langle f_{x_\ell}, \psi_j\rangle^2,
\end{equation}
since 
$$ \sum_{j=1}^{D_U}\langle f_{x_m}, \psi_j\rangle^2= \int_{-\lambda}^{\lambda} |f_{x_m}(u)|^2 du=1,$$
by Parseval's theorem and \eqref{Eq:IntegralFz}. Thus we deduce that in this case we have 
\begin{equation}\label{Eq:LowerBoundProjectionReal}
\sum_{j=1}^{D_U} \langle f_{x_{\ell}}, \psi_j\rangle^2\geq 1+O(\tau^2)>\frac12.
\end{equation}

Now suppose that $|x_{\ell}-\re(z_m)|\leq \tau $ for some $1\leq m\leq q$. Since $g_{z_m}\in U$, it follows similarly to \eqref{Eq:CSSmallDistance} that 
\begin{equation}\label{Eq:LowerBoundProjectionComplex} \sum_{j=1}^{D_U}\langle f_{x_\ell}, \psi_j\rangle^2\geq \frac{1}{\int_{-\lambda}^{\lambda} |g_{z_m}(u)|^2 du} \langle f_{x_{\ell}}, g_{z_m}\rangle^2.
\end{equation}
Write $z_m= x+iy$ for some $x, y\in \R$, and recall our assumption that $|x_{\ell}-x|\leq \tau.$ Then we have 
\begin{align}
\langle f_{x_{\ell}}, g_{z_m}\rangle&=\frac12 \Big(\langle f_{x_{\ell}}, f_{z_m}\rangle+ \langle f_{x_{\ell}}, f_{\overline{z}_m}\rangle\Big)=\int_{-\lambda}^{\lambda} \eta(u)^2 e^{-2\pi i (x_{\ell}-x)u} \cosh(2\pi y u)du \nonumber\\
&= (1+O(\tau^2))K(iy), 
\end{align}
since $$\int_{-\lambda}^{\lambda} \eta(u)^2 \cosh(2\pi y u) du= K(iy)>0,$$
which follows from the fact that $\eta$ is even.
Moreover, since $\langle f_z, f_{\overline{z}}\rangle=1$ by \eqref{Eq:KernelFactorization} we get 
$$\int_{-\lambda}^{\lambda} |g_{z_m}(u)|^2 du= \frac{1}{2}+\frac{1}{4} \left(\int_{-\lambda}^{\lambda} |f_{z_m}(u)|^2 du+ \int_{-\lambda}^{\lambda} |f_{\overline{z}_m}(u)|^2 du\right)= \frac{1}{2}+ \frac{K(2iy)}{2}.$$
Inserting these estimates in \eqref{Eq:LowerBoundProjectionComplex} implies that 
\begin{equation}\label{Eq:LowerBoundProjectionComplex2}
\sum_{j=1}^{D_U}\langle f_{x_\ell}, \psi_j\rangle^2\geq (1+O(\tau^2)) \frac{2K(iy)^2}{1+K(2iy)}\geq \frac12-  \max\left(\frac12-\frac{2K(iy)^2}{1+K(2iy)}, 0\right) +O(\tau^2),
\end{equation}
since 
\begin{align*}K(iy)^2 \leq \left(\int_{-\lambda}^{\lambda} \eta^2(u) du\right)\left(\int_{-\lambda}^{\lambda} \eta^2(u) \cosh(2\pi yu)^2du\right)= \frac{1+K(2iy)}{2},
\end{align*}
by the Cauchy--Schwarz inequality. 
Now as before we split the set $\mc{I}$ defined in \eqref{Eq:DefinitionIsolated} into a disjoint union $\mc{I}_1\cup \mc{I}_2$ where $\mc{I}_1$ is the set of those $x\in \mc{I}$ such that $|x-\re(z)|>\tau$ for all $z\in \mc{Z}\setminus{\{x\}}$. Our goal is to find a lower bound for $|\mc{I}_1|.$ If $x_{\ell}\in \mc{I}_2$ then either there is an integer $m$ such that $n+1\leq m \leq n+r$ and $|x_{\ell}-x_m|\leq \tau$ in which case we use \eqref{Eq:LowerBoundProjectionReal}, or there is an integer $m$ such that $1\leq m\leq q$ and $|x_{\ell}-\re(z_m)|\leq \tau$ in which case we use \eqref{Eq:LowerBoundProjectionComplex2}. Thus we deduce that 
$$S_0\geq \sum_{x_{\ell}\in \mc{I}_2}\sum_{j=1}^{D_U}\langle f_{x_\ell}, \psi_j\rangle^2 \geq \frac{|\mc{I}_2|}{2} -\sum_{\ell=1}^q \max\left(\frac12-\frac{2K\big(i\im(z_{\ell})\big)^2}{1+K\big(2i\im(z_{\ell})\big)}, 0\right) + O(\tau^2 N).$$
Inserting this in \eqref{Eq:MainInequalitySmallTau} and using that $d=|\mc{I}|= |\mc{I}_1|+|\mc{I}_2|$ completes the proof.

\end{proof}

\section{Proofs of Theorems \ref{Thm:Isolated0}  and \ref{Thm:Kconcecutive}}

For $T\geq 20$, let $\mc Z_T$ be the finite multiset
$$
\mc Z_T
:=
\left\{
-i\left(\rho_n-\frac12\right)\frac{\log T}{2\pi}
:
0<\gamma_n\leq T
\right\},
$$
where every zero $\rho_n$ occurs with its multiplicity. The functional
equation of the zeta function implies that $1-\overline{\rho_n}$
is a zero with the same multiplicity as $\rho_n$. Moreover, we have
$$
\overline{
-i\left(\rho_n-\frac12\right)\frac{\log T}{2\pi}
}
=
-i\left(1-\overline{\rho_n}-\frac12\right)
\frac{\log T}{2\pi}.
$$
Thus, $\mc Z_T$ is invariant under complex conjugation. Furthermore,
an element of $\mc Z_T$ is real if and only if the corresponding zero $\rho_n$
satisfies $\beta_n=1/2$. To prove Theorems \ref{Thm:Isolated0} and \ref{Thm:Kconcecutive} we will apply Propositions \ref{Pro:IsolatedZeros}, \ref{Pro:KConsecutiveGaps} and \ref{Pro:IsolatedSMALL} to the multiset $\mc{Z}_T$.  We first construct the kernel that we shall use throughout. This is an approximation of the Montgomery-Taylor kernel and follows our construction in Section 3 of \cite{La}. Let
\begin{equation}\label{Eq:Definitionf_0}
f_0(x)
:=\begin{cases}
\frac{\cos(\sqrt2x)}
{\sqrt2\sin(1/\sqrt2)} &\textup{ if } x\in [-1/2, 1/2],\\
0 & \textup{ otherwise.}
 \end{cases}
\end{equation}
Since
$
|\sqrt2x|\leq\frac1{\sqrt2}<\frac{\pi}{2}$ for all $x\in [-1/2, 1/2]$,
the function $f_0$ is strictly positive on $[-1/2, 1/2]$. Moreover, $f_0$ is real valued and even, and we have
$$
\int_{-\infty}^{\infty}f_0(x)\,dx
=
\frac{1}{\sqrt2\sin(1/\sqrt2)}
\int_{-1/2}^{1/2}\cos(\sqrt2x)\,dx
=1.
$$
Let $K_0(z)= \widehat{f_0}(z)$ for $z\in \C$.  
Let $0<\varepsilon<1/10$ be small and fixed, and take
$
\psi_{\eps}\in C_c^\infty((-1/2,1/2))
$
 be an even function such that
$
0\leq\psi_\eps\leq1
$ and
$
\psi_\eps(x)=1$
for
$|x|\leq\frac12-\eps.$
Put
\begin{equation}\label{Eq:DefAEpsilon}
A_\eps
:=
\int_{-1/2}^{1/2}\psi_\eps(x)^2f_0(x)\,dx.
\end{equation}
We now define
$$
\eta_\eps(x)
:=
\frac{\psi_\eps(x)\sqrt{f_0(x)}}{\sqrt{A_\eps}}
\quad \textup{ and } \quad 
f_\eps(x):=\eta_\eps(x)^2
=
\frac{\psi_\eps(x)^2f_0(x)}{A_\eps}.
$$
Since $f_0$ is smooth and strictly positive on $(-\frac12, \frac12)$ it follows that
$
\eta_\eps,f_\eps\in C_c^\infty\big((-1/2,1/2)\big).
$
Moreover, $\eta_\eps$ and $f_\eps$ are real valued and even, $f_\eps\geq0$ and
$
\int_{-\infty}^{\infty}f_\eps(x)\,dx=1.
$
For all $z\in \C$ we define the kernel
\begin{equation}\label{Eq:DefinitionKernel}
K_\eps(z):=\widehat{f_\eps}(z)
=\widehat{\eta_\eps^2}(z).
\end{equation}
Then $
K_\eps(0)=1
$ and for all $u\in \R$ we have 
\begin{equation}\label{Eq:ExpressionKernelEpsilon}
K_\eps(u)= \frac{2}{A_{\eps}}\int_{0}^{1/2} \psi_{\eps}(x)^2f_0(x) \cos(2\pi x u) dx,
\end{equation}
since $\psi_{\eps}$ and $f_0$ are both even functions. We will need the following lemma from \cite{La} which follows from the unconditional pair correlation estimate of Baluyot, Goldston, Suriajaya and Turnage-Butterbaugh \cite[Lemma 5]{BGST}.

\begin{lem}[Lemma 3.2 of \cite{La}]\label{Lem:ConstructionKernel}
Put
\begin{equation}\label{Eq:MontgomeryTaylorConstant}
C_{\mathrm{MT}}
:=
\frac12+\frac1{\sqrt2}
\cot\left(\frac1{\sqrt2}\right)
=
1.327499\ldots.    
\end{equation}
Let $0<\varepsilon<1/10$ be small and fixed. We have
$$
\sum_{\substack{\rho,\rho'\\0<\gamma,\gamma'\leq T}}
K_\varepsilon\left(
i(\rho-\rho')\frac{\log T}{2\pi}
\right)^2
=
\left(
C_{\eta_\varepsilon}
+
O_{\varepsilon}\left(
\frac{1}{\sqrt{\log T}}
\right)
\right)
\frac{T}{2\pi}\log T,
$$
where
$
\left|C_{\eta_\varepsilon}-C_{\mathrm{MT}}\right|
<\varepsilon.
$
\end{lem}

\subsection{Proof of the estimate \eqref{Thm:Isolated}} Let $K_0= \widehat{f_0}$  where $f_0$ is defined in \eqref{Eq:Definitionf_0}. Note that for $u\in \R$ we have 
\begin{equation}\label{Eq:KernelDecreasing}
    K_0(u)= 2\int_{0}^{1/2}f_0(x)\cos(2\pi u x) dx, \text{ and } K_0'(u)= -4\pi \int_{0}^{1/2}xf_0(x)\sin(2\pi u x) dx< 0,
\end{equation}
for all $u\in (0, 1]$. Thus $K_0$ is strictly decreasing on $[0, 1]$. Therefore, for $0\leq \tau\leq 1/2$ we have $\min\{ K_0(u): u\in [-2\tau, 2\tau]\}= K_0(2\tau).$ Let $\tau_0$ be the first positive solution to $K_0(2\tau)=5/6.$ Then one has 
\begin{equation}\label{Eq:DefinitionTau0}
\tau_0=  0.169641\ldots>\frac{1}{6}.
\end{equation}
Since $\tau_0\in [0,1/2]$ we deduce that $\kappa_0=\min\{ K_0(u): u\in [-2\tau_0, 2\tau_0]\}= K_0(2\tau_0)=5/6$.

Let $0<\eps<1/10$ be fixed and small and $K_{\eps}$ be defined in \eqref{Eq:DefinitionKernel}. Similarly to \eqref{Eq:KernelDecreasing} we have 
$$ K_\eps'(u)= -\frac{4\pi}{A_{\eps}} \int_{0}^{1/2}x\psi_{\eps}(x)^2f_0(x)\sin(2\pi u x) dx< 0$$ for all $u\in (0, 1]$. Thus, it suffices to show that $K_\eps(2\tau_0)\geq 5/6$ in order to apply Proposition \ref{Pro:IsolatedZeros} with $K=K_{\eps}$ and $\tau=\tau_0$. Since $K_0(2\tau_0)=5/6$ and $K_{0}(0)=1$  it follows from \eqref{Eq:ExpressionKernelEpsilon} that 
\begin{align}\label{Eq:OptimalityKTau_0}\frac{A_{\eps}}{2}\left(K_{\eps}(2\tau_0)-\frac56\right)&=  \int_{0}^{1/2}\psi_\eps(x)^2f_0(x)\cos(4\pi x\tau_0)dx- \frac56\int_{0}^{1/2} \psi_\eps(x)^2f_0(x) dx \nonumber\\
& = \int_{0}^{1/2}\big(\psi_\eps(x)^2-1\big)f_0(x)\cos(4\pi x\tau_0)dx +\frac{5}{12}- \frac56\int_{0}^{1/2} \psi_\eps(x)^2f_0(x) dx \nonumber\\
&= \int_{0}^{1/2}\big(1-\psi_\eps(x)^2\big)f_0(x)\left(\frac{5}{6}-\cos(4\pi x\tau_0)\right)dx,\\
&= \int_{1/2-\eps}^{1/2}\big(1-\psi_\eps(x)^2\big)f_0(x)\left(\frac{5}{6}-\cos(4\pi x\tau_0)\right)dx,\nonumber
\end{align}
since $\psi_\eps(x)=1$ for $0\leq x\leq 1/2-\eps$. Now note that $1/6<\tau_0<1/4$ and hence for all real numbers $x$ such that  $2/5<1/2-\eps\leq x\leq 1/2$ we have $\cos(4\pi x\tau_0)\leq 5/6$. Thus, we deduce that $K_{\eps}(2\tau_0)\geq 5/6$. Therefore, applying Proposition \ref{Pro:IsolatedZeros} to the multiset $\mc{Z}_T$ with $K=K_{\eps}$ and $\tau=\tau_0$ we obtain 
\begin{align*}
&\left|\left\{0<\gamma_n\leq T: \min(\gamma_{n+1}-\gamma_n, \gamma_n-\gamma_{n-1})>\frac{2\pi\tau_0}{\log T} \right\}\right|\\
&\geq \frac52 N(T) -
\frac32\sum_{\substack{\rho_n,\rho_m\\0<\gamma_n,\gamma_m\leq T}}
K_\varepsilon\left(
i(\rho_n-\rho_m)\frac{\log T}{2\pi}
\right)^2\geq \left(\frac{5}{2}-\frac{3}{2}C_{\eta_\varepsilon}+o(1)\right)N(T).\end{align*}
Since 
$
\left|C_{\eta_\varepsilon}-C_{\mathrm{MT}}\right|
<\varepsilon$ we deduce that
$$\liminf_{T\to \infty}\frac{1}{N(T)}\left|\left\{0<\gamma_n\leq T: \min(\gamma_{n+1}-\gamma_n, \gamma_n-\gamma_{n-1})>\frac{2\pi\tau_0}{\log T} \right\}\right|\ge \frac{5}{2}-\frac{3}{2}C_{\mathrm{MT}}-\frac32 \eps.$$
Letting $\eps\to 0$ and noting that $5/2-3/2C_{\mathrm{MT}}= 0.508751\dots $ completes the proof. 


\subsection{Proof of Theorem \ref{Thm:Kconcecutive}} The proof is similar to that of  \eqref{Thm:Isolated}. First, we observe that the lower bound that follows from Proposition \ref{Pro:KConsecutiveGaps} for $\mc{Z}_T$ in the case $k\geq 6$ is negative since $(k+3)/2- (k+1)C_{\mathrm{MT}}/2<0$ for $k\geq 6$. Therefore, we only consider the cases $k=3, 4, 5$. For such a $k$ let $\tau_k$ be the smallest positive solution to the equation
$$ K_0(2\tau_k)=1-\frac{1}{2(k+1)}.$$
Then we have 
$$ \tau_3= 0.145863\dots>\frac{1}{7}, \quad \tau_4= 0.129916\dots >\frac{1}{8} , \quad  \textup{ and } \quad \tau_5=0.118269\dots>\frac19.$$
Let $0<\eps<1/10$ be small and fixed. Then a similar argument leading to \eqref{Eq:OptimalityKTau_0} gives $K_{\eps}(2\tau_k)\geq 1-1/(2(k+1))$. Since $K_{\eps}$ is decreasing on $[0, 1]$, one has $K_{\eps}(x)\geq 1-1/(2(k+1))$ for all real numbers $x$ such that $|x|\leq 2\tau_k$. Therefore, applying Proposition \ref{Pro:KConsecutiveGaps} to the multiset $\mc{Z}_T$ we deduce that the number of the ordinates $0<\gamma_n\leq T$ such that 

\noindent $\min\big(\gamma_{n+k}-\gamma_{n+k-1}, \dots, \gamma_{n+1}-\gamma_n\big)> 2\pi \tau_k/\log T$ is 
\begin{align*}
&\geq \frac{k+3}{2}N(T)- \frac{k+1}{2}\sum_{\substack{\rho_n,\rho_m\\0<\gamma_n,\gamma_m\leq T}}
K_\varepsilon\left(
i(\rho_n-\rho_m)\frac{\log T}{2\pi}
\right)^2+O(1)\\
&\geq \left(\frac{k+3}{2}- \frac{k+1}{2} C_{\eta_{\eps}}+o(1)\right) N(T).
\end{align*}
Proceeding as in the proof of  \eqref{Thm:Isolated} (by taking the $\liminf_{T\to \infty}$ of the ratio by $N(T)$ and then letting $\eps\to 0$) completes the proof. 

\subsection{Proof of the estimate \eqref{Thm:IsolatedSmall}}
To deduce  \eqref{Thm:IsolatedSmall} from Proposition \ref{Pro:IsolatedSMALL} we also need to control the real parts of the off critical line zeros of the zeta function. We will achieve this using the following explicit zero density result due to Conrey \cite{CoBLMS}. 
\begin{lem}\label{Lem:Conrey}
Let $\rho=\beta+i \gamma$ denote a non-trivial zero of the zeta function. We have
$$ \sum_{\substack{0<\gamma\leq T\\ \beta>1/2}} \left(\beta-\frac12\right)\leq \big(0.0806 +o(1)\big) T.
$$
\end{lem}
\begin{proof}
This corresponds to \cite[Theorem 3]{CoBLMS} since the left hand side equals
$$ \int_{1/2}^1 N(\sigma, T) d\sigma,$$
where $N(\sigma, T)$ is the number of non-trivial zeros $\rho$ with $0<\gamma\leq T$ and $\beta>\sigma$. Although Conrey did not include a proof of this result in \cite{CoBLMS}, its proof follows exactly along the lines of the proof of the Corollary in the work of Balasubramanian, Conrey and Heath-Brown \cite{BaCoHB} (who obtained the slightly weaker constant $0.0845$) by replacing the mollifier in their result (which has length $T^{9/17-\eps}$) by Conrey's  mollifier in \cite{Co89} which has length $T^{4/7-\eps}$.
\end{proof}
\begin{proof}[Proof of \eqref{Thm:IsolatedSmall}]
Let $c>0$ be small. Let $0<\eps<1/10$ be a suitably small parameter to be chosen. Applying Proposition \ref{Pro:IsolatedSMALL} to the multiset $\mc{Z}_T$ with $\tau=c$ and $K=K_{\eps}$ gives 
\begin{align}\label{Eq:SMallGapsProof}
& \left|\left\{0<\gamma_n\leq T: \min(\gamma_{n+1}-\gamma_n, \gamma_n-\gamma_{n-1})>\frac{2\pi c}{\log T} \right\}\right|\nonumber\\
& \geq 2N(T) -\sum_{\substack{\rho_n,\rho_m\\0<\gamma_n,\gamma_m\leq T}}
K_\varepsilon\left(
i(\rho_n-\rho_m)\frac{\log T}{2\pi}
\right)^2 \\
& \qquad -\sum_{\substack{0<\gamma_n\leq T\\ \beta_n>1/2}} \max\left(1-\frac{4K_{\eps}\left(i \left(\beta_n-\frac12\right)\frac{\log T}{2\pi}\right)^2}{1+K_{\eps}\left(2i \left(\beta_n-\frac12\right)\frac{\log T}{2\pi}\right)}, 0\right) +O\left(c^2 N(T)\right), \nonumber
\end{align}
since $K_{\eps}(iy)$ is even for $y\in \R$. By Lemma \ref{Lem:ConstructionKernel}
it suffices to estimate the second sum on the right hand side of \eqref{Eq:SMallGapsProof}. Let $y\geq 0$ be a real number. Since $A_{\eps}=1+O(\eps)$ and $0\leq \psi_{\eps}(x)\leq 1$ for all $x$ with $\psi_{\eps}(x)=1$ when $0\leq x\leq 1/2-\eps$, it follows from the definitions of $K_{\eps}$ and $K_0$ that 
\begin{align}\label{Eq:ApproxKEpsilon}
K_{\eps}(iy)&= \frac{2}{A_{\eps}}\int_{0}^{1/2} \psi_{\eps}(x)^2f_0(x) \cosh(2\pi x y) dx\nonumber\\
&= 2\int_{0}^{1/2} \psi_{\eps}(x)^2f_0(x) \cosh(2\pi x y) dx +O(\eps \cosh(\pi y))\\
&= K_{0}(iy)+O(\eps \cosh(\pi y)).\nonumber 
\end{align}
Furthermore, an easy calculation gives 
$$ K_0(iy)
= 2\int_0^{1/2} f_0(x) \cosh(2\pi yx) dx=\frac{\cosh(\pi y)
+ \sqrt{2}\,\pi y \cot(1/\sqrt{2})\sinh(\pi y)}
{1+2\pi^2 y^2}.
$$
For $y\in \R$ we define 
$$ H_0(y):= 1- \frac{4K_0(iy)^2}{1+K_0(2i y)}, \quad \text{ and } \quad  H_{\eps}(y):= 1- \frac{4K_\eps(iy)^2}{1+K_\eps(2i y)}.$$
Both $H_0$ and $H_{\eps}$ are even functions of $y$ and satisfy $H_0(y)<1$ and $H_{\eps}(y)<1$ for all $y\in \R$. Moreover, we have $H_0(0)=H_\eps(0)=-1$ and $H_{\eps}(2)= H_0(2)+O(\eps)>0.4092+O(\eps)>2/5$ if $\eps$ is suitably small. Therefore we deduce that for all $y\geq 5$ we have
\begin{equation}\label{Eq:BoundMaxRegion}
\frac{H_{\eps}(y)}{y} <\frac{1}{y} \leq  \frac15 <\frac{H_{\eps}(2)}{2} \quad \textup{ and } \quad  \frac{H_{0}(y)}{y} <\frac{1}{y} \leq  \frac15 <\frac{H_{0}(2)}{2}. 
\end{equation}
Moreover, we have $\lim_{y\to 0+} H_0(y)/y=\lim_{y\to 0+} H_\eps(y)/y=-\infty$. Since $H_0$ and $H_{\eps}$ are smooth it follows that $\alpha_0:= \max_{y>0} H_0(y)/y$ and $\alpha_{\eps}:= \max_{y>0} H_{\eps}(y)/y $ exist. Therefore, using Lemma \ref{Lem:Conrey} we deduce that 
\begin{align}\label{Eq:BoundSumOffLineZeros}
\sum_{\substack{0<\gamma_n\leq T\\ \beta_n>1/2}} \max\left(1-\frac{4K_{\eps}\left(i \left(\beta_n-\frac12\right)\frac{\log T}{2\pi}\right)^2}{1+K_{\eps}\left(2i \left(\beta_n-\frac12\right)\frac{\log T}{2\pi}\right)}, 0\right) &\leq \alpha_{\eps} \frac{\log T}{2\pi}\sum_{\substack{0<\gamma_n\leq T\\ \beta_n>1/2}}\left(\beta_n-\frac12\right)\nonumber\\
& \leq \left(0.0806\alpha_{\eps}+o(1)\right)N(T).
 \end{align}
Thus it remains to obtain a numerical approximation for $\alpha_{\eps}$.
Since $\cosh(y)$ is increasing on $[0, \infty)$ we get 
$$ 1= 2\int_{0}^{1/2}f_0(x)dx\leq K_0(iy)\leq \cosh(\pi y) \left(2\int_{0}^{1/2}f_0(x)dx\right)=\cosh(\pi y),
$$
and similarly
$$ 1= \frac{2}{A_{\eps}}\int_{0}^{1/2}\psi_{\eps}(x)^2f_0(x)dx\leq K_\eps(iy)\leq \cosh(\pi y),
$$
Hence we deduce that for $y\in [0, 1/10]$ we have 
$$ \max\big(H_0(y), H_{\eps}(y)\big) \leq 1-\frac{4}{1+\cosh(2\pi y)}\leq 1-\frac{4}{1+\cosh(\pi/5)}<0.$$
Combining these bounds with \eqref{Eq:ApproxKEpsilon} and \eqref{Eq:BoundMaxRegion} and using Mathematica to compute the maximum of $H_0(y)/y$ over the interval $[1/10, 5]$ we obtain
$$ \alpha_0= \max_{y>0} \frac{H_0(y)}{y}=\max\left\{\frac{H_0(y)}{y} : y\in [1/10, 5]\right\}< 0.21462,$$
and hence
\begin{align*}
\alpha_\eps=\max\left\{\frac{H_\eps(y)}{y} : y\in [1/10, 5]\right\}=\alpha_0+O(\eps)\leq 0.21462+O(\eps).
\end{align*}
Inserting this estimate in \eqref{Eq:BoundSumOffLineZeros} and then back in \eqref{Eq:SMallGapsProof},  and using Lemma \ref{Lem:ConstructionKernel} together with the fact that $2-C_{\mathrm{MT}}- 0.0806\alpha_0+O(\eps)>0.655202$ for a suitably small $\eps$ completes the proof.

\end{proof}


\section{Zeros of linear combinations of $L$-functions}

\subsection{A refinement of the Bombieri--Hejhal method} Suppose that $F$ is a certain linear combination of the $L$-functions $L_1, \dots, L_J$ for $J\geq 2$. The main idea of the work of Bombieri and Hejhal \cite{BoHe} is to exploit the fluctuation of the values of the $L$-functions $\{L_j\}_{j\leq J}$ on the critical line, implied by Selberg's central limit theorem, to insure that for almost all $t$ one such $L$-function dominates in a short interval around $t$, and hence its sign changes will produce zeros of the linear combination $F$. However, in order to make this idea work they need to insure that the zeros of each of the $L$-functions $L_1, \dots, L_J$ are not too clustered.  

Let $q\geq 2$ be a positive integer and $\chi$ a primitive Dirichlet character modulo $q$. Put $a=\frac{1-\chi(-1)}{2}\in\{0,1\}.$ The completed $L$-function attached to $\chi$ is
$$\Lambda(s,\chi)
  :=\left(\frac{q}{\pi}\right)^{(s+a)/2}
    \Gamma\left(\frac{s+a}{2}\right)L(s,\chi).$$
It satisfies the functional equation
\begin{equation}\label{Eq:FunctionEquation}
\Lambda(s,\chi)=\varepsilon_{\chi}\Lambda(1-s,\overline{\chi}),   
\end{equation}
where $\varepsilon_\chi=\tau(\chi)/(i^a\sqrt{q})$ is the root number and $\tau(\chi)$ is the Gauss sum attached to $\chi$. 

Let $\chi_1$, $\chi_2$ be distinct primitive characters with the same parity modulo $q$. To simplify our notation we put  
$
L_1(s):=L(s, \chi_1)$,  and $L_2(s):=L(s, \chi_2)$. 
We consider the linear combination 
$$ 
F(s):= b_1 e^{i\theta_1} L_1(s)+ b_2 e^{i\theta_2}  L_2(s),
$$
where $\theta_j=- (\arg \varepsilon_{\chi_j})/2,$
and $b_1, b_2$ are non-zero real numbers.  
For $\re(s)>1$ we can write $F(s)$ as a Dirichlet series
$ F(s)= \sum_{n=1}^{\infty} a_F(n)/n^s$
 where $a_F(n)= b_1 e^{i\theta_1} \chi_1(n)+ b_2 e^{i\theta_2}\chi_2(n)$. Since the first nonzero term dominates the others, $F(s)$ has no zeros in $\re(s)>A$ for some constant $A>0$. Hence, we may define 
$$\sigma_F := \sup \{\re(\rho)  : F(\rho ) = 0 \} \leq A .$$
Moreover, it follows from  \eqref{Eq:FunctionEquation} that $F$ satisfies the functional equation 
\begin{equation}\label{Eq:FunctionEquationF}
\Lambda_F(s)=\overline{\Lambda_F(1-\overline{s})} \quad \textup{ where } \quad \Lambda_F(s)=\left(\frac{q}{\pi}\right)^{(s+a)/2}
    \Gamma\left(\frac{s+a}{2}\right) F(s).
\end{equation}
Thus, $F$ has only ``trivial'' zeros in $\re(s) < 1- \sigma_F $. All the other zeros are in the strip $ 1- \sigma_F \leq \re(s) \leq \sigma_F $ and we may call them the nontrivial zeros of $F$. The number of nontrivial zeros $ \beta+i\gamma $ of $F(s)$ with $ 0 < \gamma \leq T$ is denoted by $N_F (T)$ and it is well-known that (see \cite[Eq. (3.12)]{BoHe})
$$ 
N_F(T) \sim \frac{T}{2\pi} \log T.
$$
For $j=1,2$ let 
$ N(T, \chi_j)$ denote the number of non-trivial zeros $\rho_n^{(j)}=\beta_n^{(j)}+i\gamma_n^{(j)}$ of $L_j(s)$ with $0<\gamma_n^{(j)}\leq T$. The Riemann von-Mangoldt formula for Dirichlet $L$-functions gives
$$ N(T, \chi_j)= \frac{T}{2\pi} \log\left(\frac{qT}{2\pi}\right)-\frac{T}{2\pi}+O(\log T).
$$
For $j=1, 2$ we define the multiset $$\mc{H}_j(T):= \left\{T< \gamma_{n}^{(j)}\leq 2T : \min(\gamma^{(j)}_{n+1}- \gamma_{n}^{(j)}, \gamma_{n}^{(j)}- \gamma_{n-1}^{(j)})\leq \frac{2\pi }{(\log T)(\log\log T)^{1/4}}\right\}, 
$$
where the ordinates are counted with multiplicity. 
By refining the method of Bombieri and Hejhal \cite{BoHe} we prove the following unconditional result. 
\begin{pro}\label{Pro:BombieriHejhal}
Let $\varepsilon>0$ be small and fixed. Let $T$ be large. The number of zeros of $F(1/2+it)$ of odd multiplicity for $t\in [T, 2T]$ is 
$$ \geq (1-\varepsilon)\frac{T}{2\pi}\log T- \sum_{\gamma_{n}^{(1)}\in \mc{H}_1(T)} 1- \sum_{\gamma_{n}^{(2)}\in \mc{H}_2(T)} 1.$$
\end{pro}
To prove this result we first record the following joint Selberg central limit theorem for $\log|L_1(1/2+it)|$ and $\log|L_2(1/2+it)|$, which is a special case of \cite[Theorem B]{BoHe}. A proof of this result was also outlined by Selberg \cite{Sel89}.
\begin{lem}\label{Lem:SelbergCLT}
As $t$ varies uniformly in $[T, 2T]$, the vector
$$\left(\frac{\log|L_1(1/2+it)|}{\sqrt{(\log\log T)/2}}, \frac{\log|L_2(1/2+it)|}{\sqrt{(\log\log T)/2}}\right)$$ converges in law to the vector $(\X_1, \X_2)$ as $T\to \infty$, where $\X_1, \X_2$ are standard normal independent random variables.  
\end{lem}
\begin{proof}
 This follows from \cite[Theorem B]{BoHe} using Selberg's zero density result for Dirichlet $L$-functions (see \cite{Sel46}), which states that  uniformly for $\sigma\geq 1/2$ we have \begin{equation}\label{Eq:ZeroDensity}
    N(\sigma, T, \chi)\ll_{\chi} T^{1-c(\sigma-1/2)}\log T
\end{equation}for some positive constant $c$, where $N(\sigma, T, \chi)$ is the number of zeros of $L(s, \chi)$ in the rectangle $\sigma<\re(s)\leq 1$ and $0<\im(s)\leq T$.
\end{proof}
As a corollary we deduce the following lemma which corresponds to the second part of Lemma 8 of \cite{BoHe} and holds unconditionally for Dirichlet $L$-functions.   
\begin{lem}\label{Lem:Separation}
Let $\delta>0$ be small and fixed. There exists a set $\mc{E}\subset [T, 2T]$ of measure $O(\delta T)$ such that for all $t\in [T, 2T]\setminus \mc{E}$ one has 
\begin{equation}\label{Eq:Separation}
\Big|\log\big|L_1(1/2+it)\big|-\log\big|L_2(1/2+it)\big| \Big|\geq \delta \sqrt{\log\log T}. 
\end{equation}
\end{lem}
\begin{proof}
By Lemma \ref{Lem:SelbergCLT} it follows that $\big(\log|L_1(1/2+it)|-\log|L_2(1/2+it)|\big)/\sqrt{\log\log T}$ converges to a standard normal distribution as $T\to \infty$. Therefore, the measure of points $t\in [T, 2T]$ for which \eqref{Eq:Separation} does not hold equals $$\left(\frac{1}{\sqrt{2\pi}}\int_{-\delta}^{\delta}e^{-x^2/2} dx +o(1)\right)T\ll \delta T, $$
as desired.
\end{proof}

Let $\eps>0$ be small and fixed. Following the notation of Bombieri and Hejhal, we now divide the interval $[T, 2T]$ into smaller intervals of size $M/\log T$, where $M$ is a suitably large constant which depends on $\eps$. Write these intervals as $[y_n, y_{n+1}]$ where  
 $y_\ell= T+\ell M/\log T$, and discard the two end side intervals since their combined measure is $\leq 2M/\log T=o(1).$ 
 The next lemma corresponds to Lemma 9 of Bombieri-Hejhal \cite{BoHe}, which follows from their Lemma 7, which is unconditional for Dirichlet $L$-functions thanks to \eqref{Eq:ZeroDensity}. 

 \begin{lem}[Lemma 9 of \cite{BoHe}]\label{Lem:BHLem9}
Let $\delta>0$ be small and fixed. For more than $(1-\delta)T(\log T)/M$ positive integers $\ell\leq T(\log T)/M$, there is a set $\mc{A}_\ell\subset[y_\ell, y_{\ell+1}]$ with measure $\textup{meas}(\mc{A}_\ell)> (1-\delta)(y_{\ell+1}-y_\ell)$
such that for $j=1, 2$ and all $t\in \mc{A}_\ell$ we have
\begin{equation}\label{Eq:StabilityZeros}
N\left(t+(1-\sqrt{\delta})\frac{M}{\log T}, \chi_j\right)- N\left(t, \chi_j\right)= (1-\sqrt{\delta})\frac{M}{2\pi}+ O\left(\frac{\sqrt{\log M}}{\delta^2}\right),
\end{equation}
and 
\begin{equation}\label{Eq:BH6.10}
\frac1M\int_0^M\Bigg|\log\bigg|L_j\Big(\frac12+ it+ i\frac{u}{\log T}\Big)\bigg|- \log\bigg|L_j\Big(\frac12+ it\Big)\bigg|\Bigg| du\ll \frac{\sqrt{\log M}}{\delta^2}.
\end{equation}
 \end{lem}
\begin{proof}[Proof of Proposition \ref{Pro:BombieriHejhal}]
Following Bombieri and Hejhal, we say that a positive integer $1\leq \ell \leq T(\log T)/M$ is ``good'' if the following conditions hold
\begin{itemize}
    \item[i)] $\ell$ satisfies the conditions of Lemma \ref{Lem:BHLem9}. 
    \item[ii)]  
    We have
    $$ \left(y_\ell, y_\ell+ \frac{\sqrt{\delta}M}{\log T}\right)\cap \mc{A}_\ell\cap\mc{E}^c\neq \emptyset.$$
\end{itemize}
Let $\mc{L}$ be the set of positive integers $\ell \leq T(\log T)/M$ which are good. By Lemma \ref{Lem:BHLem9} we have 
\begin{align*}
\textup{meas}\left(\Big(y_\ell, y_\ell+ \frac{\sqrt{\delta}M}{\log T}\Big)\cap \mc{A}_\ell\right)&\geq \textup{meas}\left(y_\ell, y_\ell+ \frac{\sqrt{\delta}M}{\log T}\right)+ \textup{meas}(\mc{A}_\ell)- (y_{\ell+1}-y_\ell)\\
&\ge \frac{\sqrt{\delta}M}{2\log T},
\end{align*}
since $0<\delta\leq 1/4$. Since the measure of $\mc{E}$ is $\ll \delta T$ by Lemma \ref{Lem:Separation} it follows from Lemma \ref{Lem:BHLem9} that the number of positive integers $\ell\leq T(\log T)/M$ which are not good is $\ll \sqrt{\delta} T(\log T)/M.$  

Let $\ell\in \mc{L}$. By condition (ii) there exists a point $t_{\ell} \in \left(y_\ell, y_\ell+ \sqrt{\delta}M/\log T\right) \cap \mc{A}_\ell\cap\mc{E}^c$. Thus \eqref{Eq:StabilityZeros} and \eqref{Eq:BH6.10} hold for $t_{\ell}$. In particular, by Markov's inequality and \eqref{Eq:BH6.10} there exists a set $\mc{B}_{\ell} \subset[0, M]$ such that 
\begin{equation}\label{Eq:MeasureB_ell}
    \textup{meas}(\mc{B}_\ell)\geq M +  O\left(\frac{M\sqrt{\log M}}{\delta^2(\log\log T)^{1/3}}\right),
\end{equation} and  for $j=1,2$ and all $u\in \mc{B}_{\ell}$ we have
\begin{equation}\label{Eq:StableL_j}
\log\bigg|L_j\Big(\frac12+ it_{\ell}+ i\frac{u}{\log T}\Big)\bigg|= \log\bigg|L_j\Big(\frac12+ it_\ell\Big)\bigg| + O\left((\log\log T)^{1/3}\right).
\end{equation}
Since $|\log|L_j(1/2+it)||$ is typically of order $\sqrt{\log\log T}$ by Selberg's central limit theorem, this shows that for many positive integers $\ell$,  $|\log|L_j(1/2+it)||$ does not oscillate a lot for most points around $t_{\ell}$. 

Let $\ell\in \mc{L}$ and consider the interval $$I_{\ell}= \left[t_\ell, t_\ell+ (1-\sqrt{\delta})\frac{M}{\log T}\right] \subset (y_{\ell}, y_{\ell+1}), $$
where the inclusion follows from the fact that $y_{\ell}< t_{\ell}< y_{\ell}+ \sqrt{\delta} M/\log T$. Hence these intervals are disjoint. Let $j_{\ell}\in \{1, 2\}$ be such that 
$$ \log \bigg|L_{j_\ell}\Big(\frac12+i t_{\ell}\Big)\bigg|= \max\Bigg(\log \bigg|L_1\Big(\frac12+i t_{\ell}\Big)\bigg|, \log \bigg|L_2\left(\frac12+i t_{\ell}\right)\bigg|\Bigg).$$
For simplicity assume that $j_{\ell}=1$, the other case being similar.
Since $t_{\ell} \notin \mc{E}$, it follows from Lemma \ref{Lem:Separation} and \eqref{Eq:StableL_j} that for all $u\in \mc{B}_{\ell}$ we have 
\begin{equation}\label{Eq:DominantL}
\log \bigg|L_1\Big(\frac12+i t_{\ell}+i\frac{u}{\log T}\Big)\bigg|\geq \log \bigg|L_2\Big(\frac12+i t_{\ell}+i\frac{u}{\log T}\Big)\bigg| + \frac{\delta}{2}\sqrt{\log\log T},
\end{equation}
if $T$ is suitably large. Let $\rho_1, \dots, \rho_{k_\ell}$ be the non-trivial zeros of $L_1(s)$ whose ordinates lie inside the interval $I_\ell$. Recall that $t_{\ell}\in \mc{A}_\ell$ and hence by \eqref{Eq:StabilityZeros} we have 
\begin{equation}
k_{\ell}= (1-\sqrt{\delta})\frac{M}{2\pi} + O\left(\frac{\sqrt{\log M}}{\delta^2}\right). 
\end{equation}
However, unlike the analogous estimate  \cite[Equation (6.20)]{BoHe} which depends on GRH, not all of the zeros $\rho_j$ are assumed to lie on the critical line. Write the ordinates of these zeros as $z_1\leq z_2\leq \dots \leq z_{k_{\ell}}$. We will only keep those zeros $\rho_m$ with $2\leq m\leq k_{\ell}-1$ whose ordinates $z_m$ are $\tau$-isolated, namely $\min\big(z_{m+1}-z_m, z_m-z_{m-1}\big)> 2\pi \tau/\log T$, where $\tau= 1/(\log\log T)^{1/4}$. Note that any $\tau$-isolated zero is necessarily simple and on the critical line. Let $\rho_m=1/2+i z_m$ be such a zero. By our choice of $\tau$ and  \eqref{Eq:MeasureB_ell} there exist $u_{m}, u_{m+1}\in \mc{B}_{\ell}$ such that $t_{\ell}+u_{m}/\log T\in (z_{m-1}, z_{m})$ and $t_{\ell}+u_{m+1}/\log T\in (z_{m}, z_{m+1})$. We will now prove that these points must force a zero of odd multiplicity for $F(1/2+it)$ for $t$ in the interval $(t_\ell+ u_{m}/\log T, t_{\ell}+u_{m+1}/\log T)$. Let  $$
\vartheta(t):=\arg\left(
 \left(\frac q\pi\right)^{(1/2+a+it)/2}
 \Gamma\left(\frac{1/2+a+it}{2}\right)\right),$$
 and put
\begin{equation}\label{Eq:RealLinear}
G(t):=e^{i\vartheta(t)} F\left(\frac12+it\right)= b_1 e^{i(\vartheta(t)+\theta_1)}L_1\left(\frac12+it\right)+ b_2 e^{i(\vartheta(t)+\theta_2)}L_2\left(\frac12+it\right).
\end{equation}
Then it follows from \eqref{Eq:FunctionEquationF} that $G$ is real valued, and the real zeros of $G$ correspond to the zeros of $F$ on the critical line. Moreover, by \eqref{Eq:DominantL} we have 
$$G\left(t_{\ell}+\frac{u}{\log T}\right)=b_1e^{i\left(\vartheta\left(t_{\ell}+\frac{u}{\log T}\right)+\theta_1\right)}L_1\left(\frac12+it_\ell+ i\frac{u}{\log T}\right) \left(1+O\left(e^{-\frac{\delta}{2}\sqrt{\log\log T}}\right)\right),
$$
for all $u\in \mc{B}_{\ell}$ and hence this holds in particular for $u_m$ and $u_{m+1}$. On the other hand, $e^{i(\vartheta(t)+\theta_1)}L_1(1/2+it)$ is real valued and has a constant sign on each of the intervals $(z_{m-1}, z_m)$ and $(z_m, z_{m+1})$ and opposite signs on the two intervals, since $1/2+iz_m$ is a simple zero of $L_1(s).$ Therefore we deduce that $G(t)$ has a sign change in the interval $(t_{\ell}+u_{m}/\log T, t_{\ell}+u_{m+1}/\log T).$ Note that these intervals are disjoint for different values of $m$. Thus we conclude that the number of zeros of odd multiplicity of $F(1/2+it)$ for $t\in I_{\ell}$ is 
$$ \geq k_{\ell}- \sum_{\gamma\in \mathcal{H}_1(T)\cap I_{\ell}} 1-2\geq (1-\sqrt{\delta})\frac{M}{2\pi} -\sum_{\gamma\in \mathcal{H}_1(T)\cap I_{\ell}} 1+ O\left(\frac{\sqrt{\log M}}{\delta^2}+1\right),$$
by \eqref{Eq:StabilityZeros}. We now sum over all the good $\ell\in \mc{L}$ which implies that the number of zeros of $F(1/2+it)$ for $t\in [T, 2T]$ is 
$$ \geq \frac{T}{2\pi}\log T- \sum_{\gamma\in \mc{H}_1(T)}1-\sum_{\gamma\in \mc{H}_2(T)}1 + O\left( \frac{\sqrt{\log M}}{\delta^2 M}T\log T+ \frac{T\log T}{M}+\sqrt{\delta} T\log T\right).
$$
Choosing $\delta= M^{-1/3}$ and $M$ to be sufficiently large in terms of $\varepsilon$ completes the proof. 
\end{proof}


\subsection{Proof of Theorem \ref{Thm:ZerosLinearComb}}
To prove Theorem \ref{Thm:ZerosLinearComb} we need to extend  \eqref{Thm:IsolatedSmall} to any Dirichlet $L$-function attached to a fixed primitive character $\chi$ modulo $q$ with $q\geq 2$. The unconditional pair correlation formula of Baluyot, Goldston, Suriajaya and Turnage-Butterbaugh \cite[Lemma 5]{BGST} has been generalized to this setting by Garunk\v{s}tis and Paliulionyt\.{e} \cite{GaPa}. The following is a dyadic variant of their Lemma 3.1.
\begin{lem}[Lemma 3.1 of \cite{GaPa}]\label{Lem:GaPa}
Let $q\geq 2$ be an integer and $\chi\pmod q$ be a primitive Dirichlet character. Let $\rho_{\chi}=\beta_{\chi}+i\gamma_{\chi}$ denote the non-trivial zeros of $L(s, \chi)$. Let $w(u)=4/(4-u^2)$ and  $r\in L^{1}(\mathbb R)$ be a real-valued even function with support in $[-1,1]$, and suppose that $r$ is Lipschitz continuous at $x =0$ (that is, $
r(x)-r(0)=O(|x|)$ as $x\to 0$). Then we have
\begin{align*}
\sum_{\substack{\rho_{\chi}, \rho_{\chi}' \\ T<\gamma_\chi,\gamma_\chi'\leq 2T}} \widehat{r}\left(i(\rho_{\chi} -\rho_\chi')\frac{\log T}{2\pi}\right) w(\rho_\chi-\rho_\chi')= \frac{T}{2\pi}\log (qT)
\left(r(0) + 2\int_{0}^1 \alpha r(\alpha)\,d\alpha+o_{r, q}(1)\right).
\end{align*}
\end{lem}
Combining this result with the proof of \cite[Lemma 3.2]{La} yields the analogue of Lemma \ref{Lem:ConstructionKernel} for Dirichlet $L$-functions. 
\begin{lem}\label{Lem:ConstructionKernelDirichlet}
Let $T$ be large, $q\geq 2$ be a fixed integer and $\chi\pmod q$ be a primitive Dirichlet character. Let $0<\varepsilon<1/10$ be small and fixed and $K_{\eps}$ be the kernel defined in \eqref{Eq:DefinitionKernel}. Then we have
$$
\sum_{\substack{\rho_\chi,\rho_\chi'\\T<\gamma_\chi,\gamma_\chi'\leq 2T}}
K_\varepsilon\left(
i(\rho_\chi-\rho_\chi')\frac{\log T}{2\pi}
\right)^2
=
\big(
C_{\eta_\varepsilon}
+
o_{\eps, q}(1)
\big)
\frac{T}{2\pi}\log T,
$$
where
$
\left|C_{\eta_\varepsilon}-C_{\mathrm{MT}}\right|
<\varepsilon
$ and $C_{\mathrm{MT}}$ is defined in \eqref{Eq:MontgomeryTaylorConstant}.
\end{lem}
Finally, we need the analogue of Lemma \ref{Lem:Conrey} for a Dirichlet $L$-function attached to a fixed primitive character. This follows along the same lines of the proof of the Corollary in \cite{BaCoHB} using the following result of Wu \cite{Wu}, who computed the second mollified moment of $L(1/2+it, \chi)$ for any Dirichlet $L$-function, with a mollifier of length $T^{\theta}$ with $0<\theta<4/7.$ 
\begin{lem}\label{LemmaWu}
Let $T$ be large, $q\geq 2$ be a fixed integer and $\chi\pmod q$ be a primitive Dirichlet character. Let $0<\theta<4/7$ and put $y=T^{\theta}.$ Define 
$$B(s,\chi):=\sum_{n\leq y}\frac{\mu(n)\chi(n)}{n^s}\frac{\log(y/n)}{\log y}
$$
Then we have 
$$ \frac1T\int_T^{2T}\left|L\left(\frac12+it, \chi\right)B\left(\frac12+it, \chi\right)\right|^2dt= 1+\frac{1}{\theta}+o_q(1). $$
\end{lem}
\begin{proof}
This follows from \cite[Theorem 1.1 (B)]{Wu} and the calculations in Section 6 of \cite{Wu} with the choice of polynomials $P_1(x)=x$, $P_2(x)=P(x)=0$.
\end{proof}
Using this result, and proceeding along the same lines as the proof of the second part of the Corollary in \cite{BaCoHB} using Littlewood's lemma, we obtain the analogue of Lemma \ref{Lem:Conrey} for Dirichlet $L$-functions.
\begin{lem}\label{Lem:AnalogueConreyL}
Let $T$ be large, $q\geq 2$ be a fixed integer and $\chi\pmod q$ be a primitive Dirichlet character. We have
$$ \sum_{\substack{T<\gamma_\chi\leq 2T\\ \beta_\chi>1/2}} \left(\beta_\chi-\frac12\right)\leq \big(0.0806 +o_q(1)\big) T.
$$
\end{lem}
We are now ready to prove Theorem \ref{Thm:ZerosLinearComb}. 
\begin{proof}[Proof of Theorem \ref{Thm:ZerosLinearComb}]
By a dyadic summation, it is enough to show that the number of zeros of odd multiplicity $\rho_F=\beta_F+i\gamma_F$ of $F(s)$ with $\beta_F=1/2$ and $T<\gamma_F\leq 2T$ is $\geq 0.31 \cdot T(\log T)/(2\pi)$ for large $T$.  Let $j=1, 2$ and consider the finite multiset
$$
\mc Z_T^{(j)}
:=
\left\{
-i\left(\rho_n^{(j)}-\frac12\right)\frac{\log T}{2\pi}
:
T<\gamma_n^{(j)}\leq 2T
\right\},
$$
where every zero $\rho_n^{(j)}$ occurs with its multiplicity. Following the proof of  \eqref{Thm:IsolatedSmall} (with $c$ being suitably small) and using Lemmas \ref{Lem:ConstructionKernelDirichlet} and \ref{Lem:AnalogueConreyL} instead of Lemmas \ref{Lem:ConstructionKernel} and \ref{Lem:Conrey} respectively, we deduce that for $j=1, 2$ we have
$$|\mc{H}_j(T)|\leq (1-\delta_0+o(1)) \frac{T}{2\pi}\log T,$$
where $\delta_0=0.6552$. Inserting this estimate in Proposition \ref{Pro:BombieriHejhal} and taking $\eps$ small enough completes the proof.
\end{proof}

\end{document}